\documentclass[12pt,twoside]{amsart}
\usepackage{amssymb}
\usepackage{amscd}
\usepackage{xy}

\xyoption{arrow}
\xyoption{matrix}
\xyoption{ps}
\xyoption{dvips}

\nonstopmode

\numberwithin{equation}{section}
\newtheorem{theorem}{Theorem}[section]
\newtheorem{lemma}[theorem]{Lemma}
\newtheorem{proposition}[theorem]{Proposition}
\newtheorem{corollary}[theorem]{Corollary}

\theoremstyle{definition}
\newtheorem{definition}[theorem]{Definition}
\newtheorem{remark}[theorem]{Remark}
\newtheorem{example}[theorem]{Example}

\newcommand\adj{\operatorname{adj}}
\newcommand\codim{\operatorname{codim}}
\newcommand\coker{\operatorname{coker}}
\newcommand\Ann{\operatorname{Ann}}

\newcommand\Ass{\operatorname{Ass}}
\newcommand\Hom{\operatorname{Hom}}
\newcommand\Ext{\operatorname{Ext}}

\newcommand\rank{\operatorname{rank}}
\newcommand\rankk{\operatorname{rank}_K}
\newcommand\depth{\operatorname{depth}}

\newcommand\im{\operatorname{im}}
\newcommand\length{\operatorname{length}}

\newcommand{\HH}{H_{\mathfrak m}}

\newcommand{\ER}{\operatorname{Ext}_R}
\newcommand{\HR}{\operatorname{Hom}_R}
\newcommand{\Proj}{\operatorname{Proj}}
\newcommand{\Epi}{\operatorname{Epi}}

\newcommand{\s}{\; | \;}
\newcommand{\fall}{\mbox{for all} ~}
\newcommand{\mif}{\mbox{if} ~}

\newcommand{\ffi}{\varphi}

\newcommand{\cM}{{\mathcal M}}

\newcommand{\fm}{{\mathfrak m}}
\newcommand{\fa}{{\mathfrak a}}
\newcommand{\fb}{{\mathfrak b}}
\newcommand{\fc}{{\mathfrak c}}

\newcommand{\fp}{{\mathfrak p}}

\newcommand {\ZZ}{\mathbb{Z}}

\newcommand {\PP}{\mathbb{P}}

\newcommand{\lCM}{locally Cohen-Macaulay}
\newcommand{\CM}{Cohen-Macaulay}

\begin{document}
\title[Liaison classes of modules]{Liaison classes of modules}
\author[Uwe Nagel]{Uwe Nagel$^{*}$}




\thanks{
$^*$ Department of Mathematics, University of Kentucky, 715 Patterson
  Office Tower, Lexington, KY 40506-0027, USA; \;
e-mail: uwenagel@ms.uky.edu}  


\begin{abstract} We propose a concept of module liaison that extends Gorenstein liaison of ideals and provides an equivalence relation among  unmixed modules over a commutative Gorenstein ring.
Analyzing the resulting equivalence classes we show that several results known for Gorenstein liaison are still true in the more general case of module liaison. In particular, we construct two maps from the set of even liaison classes of modules of fixed codimension into stable equivalence classes of certain reflexive modules. As a consequence, we show that the intermediate cohomology modules and properties like being perfect, Cohen-Macaulay, Buchsbaum, or surjective-Buchsbaum are preserved in even module liaison classes. Furthermore, we prove that the module liaison class of a complete intersection of codimension one consists of precisely
 all perfect modules of codimension one.
\end{abstract}


\maketitle

\tableofcontents


\section{Introduction} \label{section-intro}

So far liaison theory can mainly be considered as an equivalence relation among equidimensional subschemes. It started with the
idea to gain information on a given curve by embedding it into a well
understood curve, a linking curve,  such that there is a residual
curve that is easier to study. The idea makes sense in any dimension
and traditionally, complete intersections were used as linking
objects. This leads to the theory of complete intersection liaison. It
has reached a very satisfactory stage for Cohen-Macaulay ideals
\cite{HU-Annals} and for  subschemes of codimension two (cf.\
\cite{rao2}, \cite{BBM}, \cite{MDP-buch}, \cite{Nollet-liai},
\cite{N-gorliaison}). 

However, it is impossible to extend all the nice results about
codimension two subschemes to higher codimension. Recently, a number
of papers (most notably \cite{KMMNP}) have shown that a more
convincing theory emerges if one allows as linking schemes instead of
complete intersections, more generally, arithmetically Gorenstein
schemes. This theory is called Gorenstein liaison. For an extensive
introduction, we refer to \cite{Migliore-buch} or \cite{torino}. 
The results in \cite{KMMNP} suggest to think of Gorenstein liaison
theory as a theory of divisors on arithmetically Cohen-Macaulay
subschemes. For example, it is shown in \cite{KMMNP} that any  two
linearly equivalent 
divisors on a smooth arithmetically Cohen-Macaulay subscheme are
Gorenstein linked in two 
steps. An application of the new theory to simplicial polytopes can be
found in \cite{MN3}. 
One can interpret this success as a consequence of enlarging the
smaller complete intersection liaison classes to the larger Gorenstein
liaison classes. 

However, despite recent efforts and many partial results (cf., e.g.,
\cite{Casanellas-Miro-Roig-1}, \cite{Casanellas-Miro-Roig-2},
\cite{Harts-co-3}, 
\cite{Harts-exp}, \cite{MN4}, \cite{NNS1}, \cite{CD-Hartsh},
\cite{C-Hartsh} \cite{Harts-div}),  Gorenstein liaison classes are not
yet well understood. In this paper, we propose to obtain a better
understanding of Gorenstein liaison and to extend the range of
applications of liaison theory 
by further enlarging Gorenstein liaison classes. To this end we
introduce a new concept of module liaison. 

There are other reasons that motivate the quest for a liaison theory
of modules. Ideals or subschemes are often studied by means of
associated modules/sheaves such as the canonical module. New insight
can be expected when modules and ideals can be treated on an equal
footing. 

Module liaison will provide a new tool for studying modules. Recently,
Casanellas, Drozd, and Hartshorne \cite{CD-Hartsh} showed that liaison
classes of codimension two ideals in a normal Gorenstein algebra $R$
are related to special maximal Cohen-Macaulay modules over
$R$. Module liaison could be helpful in investigating such modules
more directly. 

The need for a liaison theory of modules is also reflected by the fact
that so far four different proposals of module linkage (including
this one) have been developed independently \cite{Yoshino-Iso},
\cite{Martin-mod-link}, \cite{Marts-S}. However, while the other
proposals do generalize complete intersection liaison, only the
concept proposed here provides an extension of Gorenstein liaison. For
a more detailed comparison we refer to Remark \ref{rem-comparison}. 

Let us now describe the structure of the paper. In Section
\ref{section-Gorenstein-modules} we introduce the modules that will be
used to link. We require that these modules have a finite self-dual
resolution. Modules with this property are called quasi-Gorenstein
because they generalize quotients of Gorenstein rings by Gorenstein
ideals, but they are Gorenstein modules only if they are maximal
modules. We provide several classes of examples in order to illustrate
the abundance of quasi-Gorenstein modules. 

Our concept of module liaison is introduced in Section
\ref{section-module-linkage}. 
We consider unmixed modules over a local Gorenstein ring and
 graded unmixed modules over a graded Gorenstein $K$-algebra where $K$ is a
field. Throughout the paper we focus on the graded case because there
additional 
difficulties occur. Nevertheless, we show for every unmixed module
$M$, each integer $j$, and every 
quasi-Gorenstein module $C$ with the same dimension as $M$ that the modules $M,
M(j)$, and $M \oplus C$ all belong to the same even liaison class (Lemma
\ref{lem-abspalten}, Lemma \ref{lem-shift}). We also discuss several
examples and the relation to the other notions of module linkage. 
Furthermore, we describe some specializations of our module
liaison. For example, the concept of submodule liaison arises if we
restrict the class of considered modules to submodules of a given free
module $F$. In the 
special case where $F = R$ is a Gorenstein ring, submodule liaison is
the same as Gorenstein liaison of ideals. 

Then we begin our investigation of the properties of linked modules.
In Section \ref{sec-Hilb} we discuss the Hilbert polynomials of linked
modules. In particular, we show that 
$$
\deg C = \deg M + \deg N
$$
if the modules $M, N$ are directly linked by the module $C$.

In order to trace structural properties under liaison we introduce
so-called resolutions of $E$-type and $Q$-type in Section
\ref{sec-E-and-N-resolutions}. Proposition
\ref{prop-E-and-Q-resolution-under-linkage} shows how the $E$-type and
$Q$-type resolutions of directly linked modules are related. It allows
us 
to define maps $\Phi$ and $\Psi$ from the even liaison classes of
modules of fixed codimension into the set of stable equivalence
classes of certain reflexive modules (Theorem
\ref{even-liaison-class-versus-stable-euivalence-class}). 
The existence of these maps immediately produces necessary conditions
for two modules being in the same even liaison class. It remains a
major problem to decide whether these maps are injective since an
affirmative answer would give a parametrization of the even liaison
classes of modules. 

Much progress in liaison theory has been driven by the question which
properties are transferred under liaison. 
In Section \ref{sec-transfer} we use the maps $\Phi$ and $\Psi$ to
extend various results in \cite{Schenzel-liai}, 
\cite{N-gorliaison}, \cite{Marts-S}. For example, we show that the
projective dimension as well as (up to degree shift) the intermediate
local 
cohomology modules are preserved in an even module liaison class. The
same kind of preservation is true for the properties being 
Cohen-Macaulay, locally Cohen-Macaulay, Buchsbaum, and
surjective-Buchsbaum, but even in the 
whole liaison class.

The final Section \ref{sec-codim-1-liaison} is devoted to the
description of a whole module liaison class. Its main result, Theorem
\ref{thm-codim1-perfect-mod}, says that $M$ is in the liaison class
of $R/a R$ where $a \neq 0$ is any element of the Gorenstein domain
$R$ if and only if $M$ is a perfect module of codimension one. Note
that this result would follow immediately if we knew that the maps
$\Phi$ and $\Psi$ were injective. 

Our concept of module liaison could easily be extended to a non-commutative setting. The resulting theory should
certainly be investigated. We leave this for future work.
\medskip

{\bf Acknowledgments}
\smallskip

The ideas of this paper have evolved over a long time. Some of our results have been first presented in a seminar at the University of Hagen, Germany, in 1999. I would like to thank Luise Unger for the invitation  and her hospitality. 

Furthermore, I also thank the referee for many helpful comments and suggestions. 

\section{Quasi-Gorenstein modules} \label{section-Gorenstein-modules}

In this section we introduce the modules we will use for linkage.

Throughout the paper $R$ denotes a local Gorenstein ring with
maximal ideal
$\fm$ or a standard graded Gorenstein $K$-algebra over the field $[R]_0 = K$. In the
latter case $\fm = \oplus_{i>0} [R]_i$ denotes the irrelevant maximal ideal
of $R$. Usually we focus on the graded case in order to keep track of
occurring degree shifts. Ignoring degree shifts, all  definitions and
results hold analogously in the local case.

Since the ring $R$ will be  fixed we often refer to $R$-modules just
as modules. Moreover, all modules will be finitely generated unless specified otherwise.

We denote the $i$-th local cohomology module of the module $M$ by $\HH^i (M)$.
We will use two  duals of $M$, the $R$-dual $M^* := \Hom_R (M, R)$ and the
Matlis dual  $M^{\vee}$. Note that the latter is the graded module
$\Hom_K (M, K)$ if $M$ is graded.

The Hilbert function $\rankk [M]_t$ of a noetherian or
artinian graded $R$-module $M$ is denoted by $h_M(t)$. The Hilbert polynomial $p_M(t)$ is
the polynomial such that $h_M(j) = p_M(j)$ for all sufficiently large $j$. The {\em
index of regularity} of $M$ is
$$
r(M) := \inf \{i \in \ZZ \s h_M(j) = p_M (j) \quad \fall j \geq i \}.
$$
The shifted module $M(j)$, $j \in \ZZ$, has the same module structure as $M$, but
its grading is given by $[M(j)]_i := [M]_{i+j}$.

Let $M$ be an $R$-module where $n+1 = \dim R$ and $d = \dim M$. Then
$$
K_M = \Ext^{n+1-d}_R(M,R) (r(R)-1) \cong \Ext^{n+1-d}_R(M,K_R)
$$
is said to be the {\itshape canonical module} of $M$.
It is  the $R$-module representing the
functor $\HH^d(M \otimes_R \_\_{})^{\vee}$.

Recall that a perfect module is a Cohen-Macaulay $R$-module with finite
projective dimension.

\begin{definition} \label{def-Gorenstein-module}
A {\it quasi-Gorenstein $R$-module} $M$ is a finitely generated, perfect $R$-module
such that there is an integer $t$ and a (graded) isomorphism $M
\stackrel{\sim}{\longrightarrow} K_M(t)$.
\end{definition}

\begin{remark} \label{rem-Gor-module} 
(i) Following Sharp \cite{Sharp}, $M$ is a Gorenstein $R$-module if its completion $\hat{M}$ is isomorphic to a direct sum of copies of $K_{\hat{R}}$. In particular, it is a maximal $R$-module. Hence, a quasi-Gorenstein module is Gorenstein if and only if it is maximal because in this case it is simply  a finitely
generated, free $R$-module.

(ii) Let $M = R/I$ be a cyclic module. Then the following conditions are
equivalent (cf., e.g., \cite{Bruns-Herzog-Buch}, Theorem 3.3.7):
\begin{itemize}
\item[(a)] $R/I$ is a quasi-Gorenstein $R$-module.
\item[(b)] $R/I$ is a Gorenstein ring and $I$ is a perfect ideal.
\item[(c)] $I$ is a Gorenstein ideal.
\end{itemize}
\end{remark}

We denote by $M^*$ the $R$-dual $\HR(M, R)$ of an $R$-module $M$. The number
$\codim M := \dim R - \dim M$ is called the {\it codimension} of $M$.

Let $M$ be a perfect module of codimension $c$ with minimal free resolution
$$
0 \to   F_c \stackrel{\ffi_c}{\longrightarrow} F_{c-1} \to \ldots
\stackrel{\ffi_1}{\longrightarrow} F_0 \to M \to 0.
$$
We call this resolution {\it self-dual} if there is an integer $s$ such that
the dual resolution
$$
0 \to F_0^*(s) \stackrel{\ffi_1}{\longrightarrow}^* F_1^*(s) \to \ldots
\stackrel{\ffi_c^*}{\longrightarrow} F_c^*(s) \to \Ext^c_R(M, R)(s) \to 0
$$
is (as exact sequence) isomorphic to the minimal free resolution of $M$.

We denote the initial degree of a graded module $M$ by
$$
a(M) := \inf \{i \in \ZZ \s [M]_i \neq 0 \}.
$$

\begin{lemma} \label{lem-Gor-versus-selfdual-res}
Let $M$ be a perfect module. Then we have:
\begin{itemize}
\item[(a)] $M$ is a quasi-Gorenstein module if and only if its minimal free
  resolution is self-dual.
\item[(b)] If $M \cong K_M(t)$ then $t = 1- r(M) - a(M)$.
\end{itemize}
\end{lemma}

\begin{proof}
(a) If $M$ has a self-dual minimal free resolution then we have in
particular $M \cong \ER(M, R)(t)$ for some integer $t$. Thus, $M$ is a
quasi-Gorenstein module. The converse follows from the uniqueness properties
of minimal free resolutions.

(b)
The Hilbert function $h_M$ and the Hilbert polynomial $p_M$ of $M$
 can be compared by means of the following Riemann-Roch type formula
$$
h_M(j) - p_M(j) = \sum_{i=0}^d  (-1)^i \, \rankk [\HH^i(M)]_j
$$
where $d = \dim M$.
Since $M \cong K_M(t)$ is Cohen-Macaulay we obtain
$$
h_M(i) - p_M(i) = \rankk [\HH^d(M)]_i = \rankk [K_M]_{-i} = h_M(-i-t).
$$
Using the definitions of $a(M)$ and $r(M)$ we deduce $r(M) = 1 - a(M) - t$.
\end{proof}

There is an abundance of quasi-Gorenstein modules though one has to be more careful in the graded case than in the local case.

\begin{remark} \label{rem-dsum-qGor} While over a local ring the direct sum of quasi-Gorenstein modules is again quasi-Gorenstein, this is not always true for graded modules. In fact, if $C$ is a graded quasi-Gorenstein module then, for example, $C^2 \oplus C(1)$ is not quasi-Gorenstein because there is no integer $j$ such that $C^2 \oplus C(1) \cong (C^2 \oplus C(-1))(j)$.

However, $C^k$ and $C \oplus C(j)$ are always quasi-Gorenstein.
\end{remark}

There are plenty of
quasi-Gorenstein modules that are not a direct sum of proper
quasi-Gorenstein submodules.

\begin{example} \label{ex-quasi-G}
(i) Let $c \geq 3, u \geq 1$ be integers and consider a sufficiently 
  general homomorphism
  $\ffi: R(-1)^{u+c-1} \to R^u$. Then its cokernel will have the
  expected codimension $c$. Denote by $C$
the symmetric power of $\coker \ffi$ of order $\frac{c-1}{2}$. Its
  resolution is given by an Eagon-Northcott complex which is easily
  seen to be self-dual. Hence $C$
 is a quasi-Gorenstein submodule of codimension $c$.

(ii)  In \cite{Grassi} Grassi defines a strong Koszul
module as a module that has a free resolution which is analogous to
the Koszul complex.  Such a module is in particular
quasi-Gorenstein. For a specific example, take two (graded)
symmetric
homomorphisms $\ffi, \psi: F(-j) \to F$ where $F$ is a free $R$-module
of finite rank such that $\ffi \circ \psi = \psi \circ \ffi$ and $\{\det
\ffi, \det \psi\}$ is a regular sequence. Then the module $C$ with the
free resolution
$$
0 \to F (-2j) \stackrel{\begin{bmatrix} - \psi \\
\ffi \end{bmatrix}}{\longrightarrow} (F \oplus F) (-j)
\stackrel{\begin{bmatrix} \ffi & \psi
  \end{bmatrix}}{\longrightarrow} F
\to C \to 0
$$
is a quasi-Gorenstein module of codimension two.

Note, that Grassi \cite{Grassi} and B\"ohning \cite{Boehning} have
obtained  some structure theorems for quasi-Gorenstein
$R$-modules of codimension at most two that also admit a ring
structure.

(iii) Every  perfect $R$-module of codimension $c$ gives rise to a
quasi-Gorenstein modules. In fact, if $j$ is any   integer then  $M
\oplus K_M (j)$ is a quasi-Gorenstein module
because
$$
\Ext^c (M \oplus K_M (j), K_R) \cong K_M \oplus K_{K_M} (-j) \cong M
(-j) \oplus K_M = (M \oplus K_M (j)) (-j).
$$
\end{example}


\section{Module linkage: definition, examples, and specializations}
\label{section-module-linkage} 

The goal of this section is to introduce our concept of module liaison and to discuss some of its variations. Finally, we will compare it with other notions of module liaison that exist in the literature.

Let $C$ be an $R$-module. We denote by $\Epi(C)$
the set of $R$-module homomorphisms $\ffi: C \to M$ where $M$ is an
$R$-module and  $\im \ffi$ has
the same dimension as $C$.  Given a homomorphisms $\ffi \in \Epi(C)$ we want
to construct a new homomorphism $L_C(\ffi)$. Ultimately, we will see that this
construction gives a map $\Epi(C) \to \Epi(C) \cup \{0\}$.

\begin{definition} \label{def-linked-homo}
Let $C$ be a quasi-Gorenstein module of codimension $c$ and let $\ffi \in \Epi
(C)$. Let $s$ be the
integer such that $\ER^c(C, R)(s) \cong C$. Consider the exact sequence
$$
0 \to \ker \ffi \to C \to \im \ffi \to 0.
$$
It induces the long exact sequence
$$
0 \to \ER^c(\im \ffi, R)(s) \to \ER^c(C, R)(s)
\stackrel{\psi'}{\longrightarrow} \ER^c(\ker \ffi, R)(s) \to \ER^{c+1}(\im
\ffi, R)(s) \to \ldots .
$$
By assumption there is an isomorphism $\alpha: C \to \ER^c(C, R)(s)$. Thus
we obtain the homomorphism $\psi := \psi' \circ \alpha: C \to \ER^c(\ker
\ffi, R)(s)$ which we denote by $L_C(\ffi)$. (Its dependence on $\alpha$ is
not made explicit in the notation.)
\end{definition}

Note that $L_C(\ffi)$ is the zero map if $\ffi \in \Epi(C)$ is injective.

In order to analyze this construction in more detail we need two preliminary
results. The first is a version of results of Auslander and Bridger
\cite{Auslander-Bridger} and Evans and Griffith \cite{EG_buch},
respectively. It is stated as Proposition 2.5 in \cite{N-gorliaison}. We
denote the cohomological annihilator $\Ann_R \HH^i(M)$ by $\fa_i(M)$.

\begin{lemma} \label{lem-k-syzch} 
Let $R$ be a Gorenstein ring and let $M$ be a finitely
generated $R$-module. Then the following conditions are equivalent:
\begin{itemize}
\item[(a)] $M$ is a $k$-syzygy.
\item[(b)]  $\dim R/\fa_i(M) \leq i - k$ for all $i < \dim R$.
\end{itemize}
Moreover, if $k \geq 3$ then conditions $\mathrm{(a)}$ and $\mathrm{(b)}$
are equivalent
to the condition that $M$ is reflexive and $\ER^i(M^*,R) = 0$ if $1 \leq i
\leq k-2$.
\end{lemma}

Recall that there is a canonical map $M \to K_{K_M}$. It is an isomorphism if $M$ is Cohen-Macaulay, but
is neither injective nor surjective, in general.

We say that $M$ is an unmixed module if all its associated prime ideals have the
same height.

\begin{lemma} \label{lem-properties-canonical-module}
Let $M$ be an $R$-module. Then we have:
\begin{itemize}
\item[(a)] Its canonical module $K_M$ is unmixed.
\item[(b)] If $M$ is unmixed then the canonical homomorphism $M \to
  K_{K_M}$ is injective.
\end{itemize}
\end{lemma}

\begin{proof}
Claim (a) is well-known. We sketch its proof because we will use also the method  for showing (b).
Let $c$ denote the codimension of $M$. Then we choose homogeneous forms
$f_1,\ldots,f_c \in \Ann M$ such that the ideal $I := (f_1,\ldots,f_c)
\subset R$ is a complete intersection. Thus, the ring $S :=R/I$ is
Gorenstein and $M$ is a maximal $S$-module. Now, we will use the fact that
$M$ is an unmixed $R$-module if and only if $M$ is torsion-free as an
$S$-module. Indeed, this follows by comparing the cohomological
characterizations of the corresponding properties (cf., for example,
Lemma \ref{lem-k-syzch} and \cite{N-gorliaison}, Lemma 2.11).

Moreover, there is an isomorphism
$$
K_M \cong \Hom_S(M, S).
$$
It implies claim (a) because the $S$-dual of a module is a reflexive
$S$-module.

Claim (b) follows similarly. Indeed, the assumption provides that $M$ is a
torsion-free $S$-module. Thus the canonical map $M \to \Hom_S(\Hom_S(M, S),
S)$ is injective. Using the isomorphism above we are done.
\end{proof}

Now we are ready to describe properties of $L_C(\ffi)$.

\begin{proposition} \label{prop-properties-of-linked-homo}
Let $C$ be a quasi-Gorenstein
module and let $\ffi \in \Epi(C)$ be a homomorphism which is not injective.
Then we have:
\begin{itemize}
\item[(a)] There is an exact sequence
$$
0 \to K_{\im \ffi}(t) \to C \to \im L_C(\ffi) \to 0
$$
where $t = 1 - r(C) - a(C)$.
\item[(b)]  $L_C(\ffi) \in \Epi (C)$.
\item[(c)] The image $\im L_C(\ffi)$ is an unmixed $R$-module.
\item[(d)] If   $\im \ffi$ is
  unmixed then there is an isomorphism $\im L_C(L_C(\ffi)) \cong
  \im \ffi$.
\end{itemize}
\end{proposition}

\begin{proof}
Let $c$ denote the codimension of $M$.

(a) According to our assumption $\ker \ffi$ is a non-trivial submodule of the
quasi-Gorenstein module $C$. Since $C$ is an unmixed module and $\Ass (\ker \ffi)
 \subset \Ass C$ we conclude that $\dim (\ker \ffi) = \dim C$.
By the definition of $\psi = L_C(\ffi)$ we know that there is an exact
sequence
$$
0 \to \ER^c(\im \ffi, R)(s) \to C
\stackrel{\psi}{\longrightarrow} \ER^c(\ker \ffi, R)(s) \to \ER^{c+1}(\im
\ffi, R)(s) \to 0
$$
where $s = t - r(R) + 1 = r(C) - r(R) - a(C)$ because of Lemma
\ref{lem-Gor-versus-selfdual-res}.  Claim (a) follows.

According to Lemma
\ref{lem-properties-canonical-module} the canonical module $K_{\ker \ffi}$ is
an unmixed module of dimension $\dim C$.  On the other hand
we have $\dim \ER^{c+1}(\im \ffi, R) < \dim \im (\ffi) = \dim C$. Hence $\im
\psi$ is an unmixed module of  dimension $\dim C$ which proves claims (b)
and (c).

(d) We use the technique of the previous lemma. Let $I \subset
\Ann (\ker \ffi)$ be a complete intersection of codimension $c$. Put $S =
R/I$.
Since $C$ is Cohen-Macaulay the exact sequence
$$
0 \to \ker \ffi \to C \to \im \ffi \to 0
$$
 induces isomorphisms
$$
\ER^i(\ker \ffi, R) \cong \ER^{i+1}(\im \ffi, R) \quad \mbox{for all}\;
i > c.
$$
By our  assumption, $\im \ffi$ is torsion-free as $S$-module. It provides that  $\dim \ER^i(\ker \ffi, R) \leq \dim R - i - 2$ for all $i >
c$ where we use the convention that the trivial module has dimension $-
\infty$. It follows that $\ker \ffi$ is a reflexive $S$-module. Hence the
canonical map $\ker \ffi \to K_{K_{\ker \ffi}} \cong \Hom_S(\Hom_S(\ker \ffi,
S), S)$ is an isomorphism.

Now we consider the exact sequence
$$
0 \to \im \psi \to K_{\ker \ffi}(t) \to \ER^{c+1}(\im
\ffi, R)(t) \to 0.
$$
We already know that $\dim (\ER^{c+1}(\im
\ffi, R)) \leq \dim C -2$. Hence the last sequence induces an isomorphism
$$
K_{K_{\ker \ffi}}(-t) \cong K_{\im \psi}.
$$
Therefore the exact sequence 
$$
0 \to K_{\im \ffi}(t) \to K_C(t) \to \im \psi \to 0
$$
provides the exact sequence
$$
0 \to  K_{\im \psi}  \to K_{K_C}(-t)  \to  K_{K_{\im \ffi}}(-t)  \to
\ER^{c+1}(\im
\psi, R) \to 0.
$$
Using the last isomorphism above we get  the following commutative diagram
with exact rows
$$
\begin{array}{ccccccl}
0 \to & K_{K_{\ker \ffi}}  & \to &  K_{K_C} &
\stackrel{\gamma}{\longrightarrow} & K_{K_{\im \ffi}} & \to \ER^{c+1}(\im
\psi, R)(t) \to 0 \\ & \uparrow & & \uparrow & & \uparrow & \\
0 \to & \ker \ffi & \to & C & \to & \im \ffi & \to 0
\end{array}
$$
where the vertical maps are the corresponding canonical
homomorphisms. Since the two leftmost vertical maps are
isomorphisms and the third one is injective we conclude that there is an
isomorphism
$$
\im \gamma \cong \im \ffi.
$$
But $\im \gamma$ is isomorphic to $\im (L_C(L_C(\ffi))$ which proves claim
(d).
\end{proof}

The preceding result allows us to define.

\begin{definition} \label{def-linking-map} 
Let $C$ be a quasi-Gorenstein module with the isomorphism $\alpha: C
\stackrel{\sim}{\longrightarrow} K_C(t)$. Then the map $L_C: \Epi(C) \to
\Epi(C) \cup\{0\}, \ffi \mapsto L_C(\ffi),$ is called the {\it linking map}
with respect to $C$ (and $\alpha$). Here $0$ denotes the trivial homomorphism
$C \to 0_R$.
\end{definition}

\begin{remark}
The linking map respects isomorphisms in the following sense:
Let $\ffi, \ffi'$ be two homomorphisms in $\Epi (C)$. Following the
description of the linking map it is not difficult to see that $\im \ffi
\cong \im \ffi'$ implies $\im L_C(\ffi) \cong \im L_C(\ffi')$. Moreover, part
(d) of the previous result shows that the converse is true provided $\im
\ffi$ and $\im \ffi'$ are unmixed modules.
\end{remark}

\begin{definition} \label{def-module-linkage}
We say that two $R$-modules $M, N$ are {\it module linked in one step} by the quasi-Gorenstein
module
$C$ if there are homomorphisms $\ffi, \psi \in \Epi(C)$ such that
\begin{itemize}
\item[(i)] $M = \im \ffi,\; N = \im \psi$ \quad and
\item[(ii)] $M \cong \im L_C(\psi),\; N \cong \im L_C(\ffi)$.
\end{itemize}
Most of the time we abbreviate module linkage by m-linkage.
\end{definition}

\begin{remark} \label{rem-shift-inv} 
(i) Modules that are module linked in one step will also be called {\em directly m-linked} modules.

(ii)
If $M$ and $N$ are directly m-linked by $C$ then, by the definition of the
linking map $L_C$ and the previous lemma, we have  that $\dim M = \dim N = \dim C$
and $M, N$ are unmixed $R$-modules.

(iii)
Module linkage is shift invariant in the following sense. The modules $M, N$ are
directly
linked by $C$ if and only if the modules $M(j), N(j)$ are directly linked by $C(j)$
where $j$ is any integer.
\end{remark}

The following observation allows us to construct plenty of modules that are linked to a given module.

\begin{remark} \label{rem-plenty-links}
(i)
If $M$ is an $R$-module such that there is an epimorphism $\ffi: C \to
M$ and $M$ has the same dimension as $C$ but is not isomorphic to $C$, then
the modules $N := \im L_C(\ffi)$ and $M$ are m-linked because $M \cong \im L_C (L_C (\ffi))$ by Proposition \ref{prop-properties-of-linked-homo}.

By abuse of notation we
sometimes write  $L_C(M)$ instead of $\im L_C(\ffi)$. Then two modules
$M$ and $N$ are m-linked in one step if and only if there is a
suitable quasi-Gorenstein module $C$ such that $N \cong L_C(M)$, or
equivalently $M \cong L_C(N)$.

(ii)
A simple way to produce an epimorphism $\ffi$ as above in order to link a given module $M$ is the following. Choose  a free $R$-module $F$ such that there is an epimorphism $\psi: F \to M$ and take a complete intersection ideal $\fc$ of codimension $c = \codim M$ in $\Ann_R M$. Then $\psi$ induces an epimorphism $F/ \fc F \to M$. Thus, the map $\ffi:  F/ \fc F \oplus K_{F/ \fc F} \to M$ where $K_{F/ \fc F}$ maps onto zero satisfies the requirements because $F/ \fc F \oplus K_{F/ \fc F}$ is quasi-Gorenstein by Example \ref{ex-quasi-G}(iii). The last step, i.e.\ adding the canonical module,  can be omitted if $F/\fc F$ is already a quasi-Gorenstein module.
\end{remark}

The following examples illustrate the flexibility of our concept of module liaison.

\begin{example} \label{ex-mlink}
(i)
Every perfect module $M$ is linked to itself as a consequence of the exact sequence
$$
0 \to K_M \to M \oplus K_M \to M \to 0.
$$
This is very much in contrast to the situation of linkage of ideals where self-linked ideals are rather rare.

(ii) Every free module $F$ of rank $r > 1$ is directly m-linked to a free module of smaller rank.

In fact, write $F = R(j) \oplus G$ and set $C := R(j) \oplus G \oplus G^*(-2j)$. Then $C^* \cong C(-2j)$, thus the exact sequence
$$
0 \to G^*(2j)  \to C \to F \to 0
$$ shows that $F = R(j) \oplus G$ is directly m-linked to $G$.
\end{example}

Allowing one more link, we can extend the  last example to non-free modules.

\begin{lemma} \label{lem-abspalten}
Let $D$ be a quasi-Gorenstein module and let $M$ be any unmixed module such that $M$ and $D$ have the same dimension. Then $M \oplus D$ can be linked to $M$ in two steps.
\end{lemma}

\begin{proof}
By assumption on $D$, there is an integer $s$ such that $K_D \cong D (-s)$.

As in Remark \ref{rem-plenty-links}, we choose a free $R$-module $F$ and a complete intersection ideal $\fc$ such that $M \in \Epi (F/\fc F)$. Then Example \ref{ex-quasi-G}(iii) shows that $C := F/\fc F \oplus K_{F/\fc F}(s)$ is quasi-Gorenstein. Thus, we can use this module to link $M$ to a module $N$.

By our choice of the twist $s$ in the definition of $C$, the module $D \oplus C$ is quasi-Gorenstein, too. It follows that the modules $D \oplus M$ and $N$ are linked by $D \oplus C$ proving our claim.
\end{proof}

By its definition, module linkage is symmetric. Thus, it generates an equivalence relation.

\begin{definition} \label{def-m-liaison} {\em Module liaison} or simply {\em liaison} is the
 equivalence relation generated by direct module linkage. Its equivalence classes are called {\em (module) liaison classes} or {\em m-liaison classes}. Thus, two modules $M$ and $M'$ belong to the same liaison class if there are modules $N_0 = M, N_1,\ldots, N_{s-1}, N_s = M'$ such that $N_i$ and $N_{i+1}$ are directly linked for all $i = 0,\ldots,s-1$. In this case, we say that $M$ and $M'$ are {\em linked in $s$ steps}. If $s$ is even then $M$ and $M'$ are said to be {\em evenly linked}.

Even linkage also generates an equivalence relation. Its equivalence classes are called {\em even (module) liaison classes}.
\end{definition}

\begin{example} \label{ex-li-class-R}
Since $R$ is linked to itself by $R^2$, the even liaison class and the liaison class of $R$ agree. It contains all non-trivial free $R$-modules of finite rank. Indeed, Example \ref{ex-mlink} show that every free module is in the liaison class of a free module of rank one. But the modules $R(j)$ ($j \in \ZZ$)  and $R$ are directly linked by $R \oplus R(j)$.

We will see in Corollary \ref{cor-free} that the finitely generated, free $R$-modules form the whole liaison class of $R$.
\end{example}

Since the module structure of a module is not changed by shifting, the following
property of module liaison is certainly desirable.

\begin{lemma} \label{lem-shift}
Let $M$ be an unmixed module and let $j$ be any integer. Then $M$ and $M(j)$ belong
to the same even m-liaison class.
\end{lemma}

\begin{proof}
Let $N$ be any module that is directly linked to $M$ by the quasi-Gorenstein module
$C$. Such modules exist by Remark \ref{rem-plenty-links}. Assume that $K_C \cong
C(-t)$. Then $C \oplus K_C (i)$ is quasi-Gorenstein for all $i \in \ZZ$ and we have the
 exact sequence
$$
0 \to K_M(t) \oplus K_C(i) \to C \oplus K_C (i) \to N \to 0.
$$
Hence, $N$ is directly linked to $M (i-t) \oplus C$ for every $i \in \ZZ$. According
to Lemma \ref{lem-abspalten}, the modules $M(i-t) \oplus C$ and $M (i-t)$ are evenly
linked. Thus, choosing $i$ appropriately we get that $N$ and $M$ as well as $N$ and
$M(j)$ are linked in an odd number of steps. Our claim follows.
\end{proof}

The following construction is most commonly used for maximal modules. We keep its name in the general case, too.

\begin{definition} \label{def-Ausl-dual}
Let $M$ be a non-free $R$-module. Let $F$ be a free $R$-module and let $\pi: F \to M$ be a minimal epimorphism, i.e.\ an epimorphism that  satisfies $\ker \pi \subset \fm \cdot F$. Then we call the  module
$$
M^{\times} := \coker \Hom_R (\pi, R)
$$
the {\em Auslander dual} of $M$. It is uniquely defined up to isomorphism.
\end{definition}

This concept will be crucial in Section \ref{sec-transfer}. Here, we just show that the Auslander dual $M^{\times}$ and $M$ belong to the same module liaison class if $M$ is a maximal module.

\begin{lemma} \label{lem-Ausl-dual-liai}
Let $M$ be a non-free, unmixed maximal $R$-module. Then  $M^{\times}$ is in the m-liaison class of $M$. 
More precisely, $M$ can be linked to $M^{\times}$ in an odd number of steps. 
\end{lemma}

\begin{proof}
Consider the following exact commutative diagram
\begin{equation*}
\begin{CD}
0 @>>> \ker \ffi \oplus F^*  @>>> F  \oplus F^* @>{\ffi}>> M  @>>> 0 \\
&& @VVV @VV{\gamma}V @VV{=}V \\
0 @>>> \ker \pi @>>> F  @>{\pi}>> M  @>>> 0 \\
\end{CD}
\end{equation*}
where $\pi$ is a minimal epimorphism, $F$ is free, and $\gamma$ is the canonical projection. Put $C := F \oplus F^*$. Then dualizing with respect to $R$ provides the exact commutative diagram
\begin{equation*}
\begin{CD}
0 @>>> M^*  @>>> F^*  @>>> M^{\times}  @>>> 0 \\
&& @VV{=}V @VV{\gamma^*}V @VVV \\
0 @>>> M^*  @>>> F \oplus F^*  @>>> \im L_C (\ffi)  @>>> 0. \\
\end{CD}
\end{equation*}
Now, the Snake lemma shows that $M$ is directly m-linked to $\im L_C (\ffi) \cong M^{\times} \oplus F$. Applying Lemma \ref{lem-abspalten} successively
we see that $M^{\times} \oplus F$ and $M^{\times}$ are evenly linked. This completes the argument.
\end{proof}

\begin{remark} \label{rem-Ausl-d}
Note that in the local case we could simply use $F$ as linking module. This shows that then $M$ and $M^{\times}$ are even directly m-linked.
\end{remark}

Before comparing our concept of module liaison with other versions of module liaison in the literature, we want to discuss some variations of our concept (cf.\ also Remark \ref{rem-comparison}).

For example, one could restrict the class of modules that are used for linkage. This would lead to (potentially) smaller liaison classes. While the definition above is designed to generalize Gorenstein liaison of ideals, allowing as linking modules only strong Koszul modules might lead to a concept of module liaison which could be viewed as the proper generalization of complete intersection liaison of ideals. We do not pursue this here.

Another variation that seems worth mentioning is to restrict the focus to submodules of a given free module.

\begin{definition} \label{def-sm-link}\  
Let $F$ be a free $R$-module. Then submodules $M', N'$ of $M$ are said
to be submodule linked, or  {\em sm-linked} for short, by the
submodule $C' \subset F$ if $F/M'$ and $F/N'$ are linked by $F/C'$. As
above, this leads to equivalence classes of unmixed submodules of
$F$. 
\end{definition}

In the very special case $F = R$, submodule liaison is equivalent to Gorenstein liaison of ideals.

\begin{lemma} \label{lem-m-vers-G-link} Two ideals $I, J$ of $R$ are sm-linked by the ideal $\fc \subset R$ if and only if $\fc$ is a Gorenstein ideal of $R$ and
$$
\fc : I = J \quad \mbox{and} \quad \fc : J = I,
$$
in other words, $I$ and $J$ are Gorenstein linked by $\fc$.
\end{lemma}

\begin{proof}
If $I, J$ of $R$ are sm-linked by the ideal $\fc$ then we have by
Proposition \ref{prop-properties-of-linked-homo}(a) the exact sequence
$$
0 \to K_{R/J} (1 - r(R/\fc)) \to R/\fc \to R/I \to 0.
$$
Thus, the isomorphism $K_{R/J} (1 - r(R/\fc)) \cong \fc : J/\fc$ shows $\fc : J = I$. Similarly, we get $\fc : I = J$, thus $I$ and $J$ are Gorenstein linked. The reverse implication is clear.
\end{proof}

In spite of the last observation we view module and submodule liaison as extensions of Gorenstein liaison of ideals.

\begin{remark} \label{rem-comparison}
There are several concepts of module liaison in the literature that have been developed independently.

The first published proposal is due to Yoshino and Isogawa
\cite{Yoshino-Iso}. They work over a local Gorenstein ring and
consider Cohen-Macaulay modules only. They say that the modules $M$
and $N$ are linked if there is a complete intersection ideal $\fc$
contained in $\Ann_R M \cap \Ann_R N$ such that $M$ is isomorphic to
the Auslander dual of $N$ considered as $R/\fc$-module. Note that we
have rephrased their definition  in a way that it makes sense also for
non-Cohen-Macaulay modules. 

Martsinkovsky and Strooker \cite{Marts-S} work in greater generality
though their main results are for modules over a local Gorenstein
ring. In this case, their definition of linkage is the same as the one
of Yoshino and Isogawa as given above. Note that this is a very
special case of our concept of linkage because the modules $M$ and $N$
are linked in the sense of the two papers mentioned above if and only
if they are m-linked by $F/\fc$ in the sense of our Definition
\ref{def-module-linkage} where $F$ is the free module in a minimal
epimorphism $F \to M$ and $\fc$ is as above some complete intersection
ideal contained in $\Ann_R M \cap \Ann_R N$. In other words, we get
the liaison concept of Martsinkovsky and Strooker by restricting
drastically the modules we allow as linking modules. Though this leads
to an extension of the concept of complete intersection liaison of
ideals it does not extend Gorenstein liaison of ideals. Another
consequence of this restriction is that the resulting liaison class of
a cyclic module $R/I$ contains only cyclic modules, thus it is
essentially just the complete intersection liaison class of $I$ when
we identify a cyclic module with its annihilator. 
 
Martin's approach \cite{Martin-mod-link} is very different. He uses
generic modules in order to link making it difficult to find any
module at all that is linked to a given one. This seems rather the opposite
of the wish for large equivalence classes. 
\end{remark}
 
In \cite{CD-Hartsh}, Hartshorne, Casanellas, and Drozd consider an 
extension of Gorenstein liaison of ideals that is not yet fully generalized
by Definition \ref{def-m-liaison}. Indeed, let $I \subset J$ be  homogenous
ideals in the polynomial ring $R = K[x_0,\ldots,x_n]$.  Then, they
define the {\em G-liaison class of $J$ in $\Proj (R/I)$} as the set of
ideals in $R$ that are sm-linked to $J$ (in the sense of Definition
\ref{def-sm-link}) such that all the ideals involved in the various
links contain $I$. If $A := R/I$ is Gorenstein we can also consider
the sm-liaison class of ideals in $A$ that is generated by $J/I$. Identifying
every ideal $\fa \subset R$ in the G-liaison class of $J$ in $\Proj
(R/I)$ with $\fa/I \subset A$, this G-liaison class is larger than the
sm-liaison class of $J$ consisting of ideals in $A$. The reason is
that, if the ideals $\fa, \fb \subset R$ are sm-linked in $R$ by $\fc$
where $I \subset \fc$, then
$\fa/I, \fb/I$ are not sm-linked in $A$ by $\fc/I$ unless $\fc/I$ has
finite projective dimension as $A$-module. This motivates the
following extension of the concepts above.

\begin{definition} \label{def-rel-liai} 
Let $A$ 
be any graded quotient ring of $R = K[x_0,\ldots,x_n]$, say $A :=
R/I$. Let $M$ be a graded 
$R$-module that is annihilated by $I$. Then we say that the $R$-module
$N$ is in the {\em  m-liaison class of $M$ relative to $I$} if $M$ 
can be 
linked to $N$ by using quasi-Gorenstein $R$-modules $C_1,\ldots,C_s$
that are all annihilated by $I$. 
\end{definition} 

 If $J \subset R$ is an ideal that
contains $I$, then, identifying an  cyclic $R$-module
with its annihilator, the m-liaison class of $R/J$ relative to
  $I$ contains the G-liaison class of $J$ in $\Proj
(R/I)$. In this sense, m-liaison relative to $I$ generalizes G-liaison in $\Proj (R/I)$. 

Furthermore, if $R/I$
is Gorenstein, then it is not too difficult to see that the m-liaison
class of $M$ relative to $I$ also 
contains the m-liaison class of $A$-modules generated by $M$ in the
sense of Definition \ref{def-m-liaison}. 

Though it seems very interesting to investigate these relative m-liaison
classes, we leave this for future work and focus on studying m-liaison
classes (cf.\ Defintion \ref{def-m-liaison}) in this paper.


\section{Hilbert polynomials under liaison}
\label{sec-Hilb}

In this section we begin to relate the properties of linked modules.
The starting point  is the following result which follows immediately
by Proposition \ref{prop-properties-of-linked-homo} (a).

\begin{lemma} \label{lem-st-seq}
If the modules $M$ and $N$ are directly  m-linked by the
quasi-Gorenstein module $C$ then there is an exact sequence of
$R$-modules
$$
0 \to K_M (t) \to C \to N \to 0
$$
where $t = 1 - r(C) - a(C)$.
\end{lemma}

As in the case of linked  ideals, there is a relation among the associated prime ideals of linked modules.

\begin{corollary} \label{cor-ass-primes}
If the modules $M$ and $N$ are directly  m-linked by  $C$ then we have
$$
\Ass_R M \cup \Ass_R N = \Ass_R C.
$$
\end{corollary}

\begin{proof}
Since linkage is symmetric we have the two exact sequences
$$
0 \to K_M (t) \to C \to N \to 0
$$
and
$$
0 \to K_N (t) \to C \to M \to 0.
$$
The claim follows because the associated primes of an unmixed module and its canonical module agree. \end{proof}

Lemma \ref{lem-st-seq} allows us  to compare the Hilbert polynomials of
linked modules.

Let $M$ be a module of dimension $d$.
If $d  > 0$ then its Hilbert polynomial can be written in the form
$$
p_M(j) = h_0(M) \binom{j}{d-1} + h_1(M) \binom{j}{d-2} + \ldots +
h_{d-1}(M)
$$
where  $h_0(M),\ldots,h_{d-1}(M)$ are integers and $h_0(M) >
0$ is called the degree of $M$.  If $\dim M = 0$ then we set $\deg M
:= \length (M)$.
By abuse of notation, the degree of
an ideal $I$ is $\deg I = h_0(R/I)$. It is just the degree of the
subscheme $\Proj (R/I)$.
Now we can state.

\begin{proposition} \label{prop-Hilbert-polynomial-and-liaison}
Let $M, N$ be graded $R$-modules  that are directly
linked by $C$. Put $s := r(C) + a(C) - 1$ and $d := \dim M$.
Then we have
\begin{itemize}
\item[(a)] $$
\deg N = \deg C - \deg M,
$$
and if in addition $d \geq 2$  then
$$
h_1(N) = \frac{s-d+2}{2} \left [\deg M - \deg N \right ] + h_1(M).
$$
\item[(b)] If  $M$ is \lCM\ then
$$
p_{N} (j) = p_{C}(j) + (-1)^{d} p_{M}(s-j).
$$
\item[(c)] If $M$ is \CM\ then
$$
h_{N} (j) = h_{C}(j) + (-1)^{d-1} [h_{M}(s-j) - p_M(s-j)].
$$
\end{itemize}
\end{proposition}

For the proof we need a cohomological characterization of the property
being unmixed.

\begin{lemma} \label{lem-char-unmixed}
The $R$-module $M$ is unmixed if and only if
$$
\dim R/\Ann_R (\HH^i(M)) < i \quad \mbox{for all} \quad i < \dim M
$$
where we define the dimension of the zero module to be $- \infty$.
\end{lemma}

\begin{proof}
Let $\{f_1,\ldots,f_d\}$ be a regular $R$-sequence in the annihilator
of $M$ where $d := \dim M$. Then the claim follows by local duality
and considering $M$ as module over $R/(f_1,\ldots,f_d)$ as in the
proof of Lemma 4 in \cite{Nagel-liaison}.
\end{proof}

Now we are ready for the proof of the proposition above.

\begin{proof}[Proof of Proposition \ref{prop-Hilbert-polynomial-and-liaison}]
Again, we use the Riemann-Roch type formula
$$
h_M(j) - p_M(j) = \sum_{i=0}^{d}  (-1)^i \, \rankk [\HH^i(M)]_j.
$$
Furthermore, we have by local duality
$$
\rankk [\HH^{d}({R/I})]_j = \rankk [K_{M}]_{-j}.
$$
Now, we show claim (c). If $M$ is Cohen-Macaulay then the formulas
above and Lemma \ref{lem-st-seq} provide
\begin{eqnarray*}
h_N(j) & = & h_C (j) - \rankk [\HH^i(M)]_{s-j} \\
& = & h_C (j) + (-1)^{d-1}[h_M (s-j) - p_M (s-j)].
\end{eqnarray*}

Having shown (c) we may and will assume for the remainder of the proof
that $d = \dim M \geq 2$. Next, we show claim (a).
According to Lemma \ref{lem-char-unmixed}, the degree of the Hilbert polynomial
of $\HH^i(M)$ is at most $\max \{0, i-2\}$. Thus, using the formulas
above  we obtain for all
$j \ll 0$
$$
\begin{array}{rcl}
-p_{M} (j) & = & (-1)^{d} \rankk [\HH^{d}({M})]_j +
o(j^{d-2}) \\[1ex]
& = & (-1)^{d} \rankk [K_{M}]_{-j} + o(j^{d-2}).
\end{array}
$$
Combined with Lemma \ref{lem-st-seq} this provides
$$
p_{N} (j) = p_{C} (j) + (-1)^{d} p_{M} (s-j) + o(j^{d-2}).
$$
Comparing coefficients we get by a routine computation
$$
\deg N = \deg C - \deg M,
$$
as claimed, and
$$
h_1(N) = (s-d+2) \deg M + h_1(M) + h_1(C). \leqno(*)
$$
Since linkage is symmetric there is an analogous formula with $M$ and
$N$ interchanged. Adding both equations provides
$$
h_1(C) = - \frac{s-d+2}{2} \deg C.
$$
Plugging this into $(*)$ we get the second statement of claim (a).

If $M$ is locally Cohen-Macaulay then $[\HH^i(M)]_j = 0$ if $i <
d$ and $j \ll 0$. Thus, an analogous (but easier) argument shows claim
(b).
\end{proof}

\begin{remark}  (i) Proposition
  \ref{prop-Hilbert-polynomial-and-liaison} generalizes Corollary 3.6
  in   \cite{N-gorliaison}.

(ii) Let us illustrate the result by considering a well-known special
 case. Consider two curves  $C_1 = \Proj (R/I)$ and $C_2 =
\Proj (R/J)$ in $\PP^n$  that are linked by a complete
intersection cut out by hypersurfaces of degree $d_1,\ldots,d_{n-1}$.
 Let us denote the arithmetic genus of the curves by $g_1$ and
$g_2$, respectively. For the linking module $C$ we have $r(C) = d_1 +
 \ldots d_{n-1} -n$ (cf., e.g., \cite{N-gorliaison}, Lemma 2.3). Thus,
 in this case Proposition
\ref{prop-Hilbert-polynomial-and-liaison}(a)  takes the familiar form
(cf. \cite{Migliore-buch}, Corollary 4.2.11)
$$
g_1 - g_2 = \frac{1}{2} (d_1 + \ldots d_{n-1} -n - 1) [\deg C_1 - \deg
C_2].
$$
\end{remark}

The next observation shows that it is easier to compare the
Hilbert functions of modules that are linked in two steps and not just
one. We will discuss more results along this line later on.

\begin{lemma} \label{Hilbert-function-under-even-liaison} Suppose $M,
  N, M'$  are graded modules  such that $M$ and $N$ are linked by
  $C$ and $N$ and $M'$ are linked by $C'$. Put $s := r(C) - r(C') +
  a(C) - a(C')$. Then we have  for all
  integers $j$:
$$
h_{M'}(j) = h_{M}(j + s) + h_{C'}(j) -
h_{C}(j + s).
$$
\end{lemma}

\begin{proof} According to Lemma \ref{lem-st-seq} we have the
  following exact sequences:
\begin{eqnarray*}
& 0 \to K_{N}(1- r(C) - a(C)) \to C \to M \to 0 & \\
& 0 \to K_{N}(1- r(C') - a(C')) \to C' \to M' \to 0. &
\end{eqnarray*}
The claim follows.
\end{proof}

In order to compare other properties and, in particular, the cohomology
of linked modules we need more tools. These will be developed in the
following section.


\section{Resolutions of $E$-type and $Q$-type}
\label{sec-E-and-N-resolutions}

The purpose of this section is to show the existence of  maps $\Phi$ and
$\Psi$  from the set
of even liaison classes into the set of stable equivalence classes of
certain reflexive modules. This will be achieved by exploiting  resolutions
of $E$-type and $Q$-type. These resolutions generalize the resolutions
of $E$-type and $N$-type of ideals (cf.\ Remark \ref{rem-to-def-of-E-and-Q-resolution} below) which have been introduced in \cite{MDP-buch}.

\begin{definition} \label{def-E-and-Q-type} 
Let $M$ be an $R$-module
   of codimension $c > 0$. Then an {\it $E$-type
    resolution}
  of $M$ is  an exact sequence of finitely generated graded $R$-modules
 $$
0 \to E \to F_{c-1} \to \ldots \to F_0 \to M \to 0
$$
where the modules $F_0, \ldots, F_{c-1}$ are free.

A {\it $Q$-type resolution} of $M$ is an exact sequence of finitely
generated graded $R$-modules
$$
0 \to G_c \to \ldots \to G_2 \to Q \to G_0 \to M \to 0
$$
where $G_0, G_2, \ldots, G_c$ are free and $\HH^i(Q) = 0$ for all $i$ with
$n+2-c
\leq i \leq n$. (Note, that for a module of codimension one a $Q$-type
resolution is the same as an $E$-type resolution.)

These resolutions of $M$ are said to be {\it minimal} if it is not possible
to split off free direct summands from any of the occurring modules
besides $M$.
\end{definition}

\begin{remark}  \label{rem-to-def-of-E-and-Q-resolution}
 A (minimal) $E$-type resolution of $M$ always exists because it is just the
beginning of a (minimal) free resolution of $M$. Thus, a minimal $E$-type
 resolution is uniquely determined up to isomorphism of
 complexes. Moreover, it follows that
$$
\HH^i(E) \cong \HH^{i-c}(M) \quad \mif i \leq n.
$$
\end{remark}

It requires some more work to show that $Q$-type resolutions exist.

\begin{lemma} \label{lem-Q-type-exists}
Every module $M$ of positive codimension admits a minimal $Q$-type resolution
$$
0 \to G_c \to \ldots \to G_2 \to Q \to G_0 \to M \to 0.
$$
It is uniquely determined up to isomorphism of complexes. Furthermore,
we have
$$
\HH^i(Q) \cong \left \{ \begin{array}{ll}
\HH^{i-1}(M) & \mif i \leq n+1-c \\
0 & \mif n+2-c \leq i \leq n.
\end{array} \right.
$$
\end{lemma}

\begin{proof}
We may assume that the codimension $c$ of $M$ is at least two.
Let
$$
G_1 \stackrel{\ffi}{\longrightarrow} G_0 \to M \to 0
$$
be a minimal presentation of $M$. Set $T := \ker \ffi$. Now consider a
so-called minimal $(c-1)$-presentation of $T$, i.e.\ an exact sequence
of graded $R$-modules
$$
0 \to P \to Q \to T \to 0
$$
such that $P$ has projective dimension $\leq c-2$,
$$
\HH^i(Q) = 0 \quad \mbox{for all $i$ with}\; n+2 - c \leq i \leq n,
$$
and it is not possible to split off a non-trivial free $R$-module
being a direct summand of $P$ and $Q$. Such a sequence exists and is
uniquely determined  by
\cite{N-loc-res}, Theorem 3.4 (cf.\ also \cite{EG_buch} in the local
case).  Using \cite{N-gorliaison}, Lemma 2.9 we see that
$$
\HH^i(Q) \cong \left \{ \begin{array}{ll}
\HH^{i-1}(M) & \mif i \leq n+1-c \\
0 & \mif n+2-c \leq i \leq n,
\end{array} \right.
$$
as claimed, and that $P$ has projective dimension $c-2$ because 
$$
\HH^{n+3-c} (P) \cong \HH^{n+2-c} (T) \cong \HH^{n+1-c} (M) \neq 0   
$$ 
if  $c \geq 3$. Hence
replacing $P$ in the exact sequence
$$
0 \to P \to Q \to G_0 \to M \to 0
$$
by its minimal free resolution provides a minimal $Q$-type resolution
of $M$.

Conversely, any $Q$-type resolution gives
rise to a $(c-1)$-presentation of $T$. Thus, the  uniqueness of the
minimal $Q$-type resolution follows from the uniqueness of the minimal
$(c-1)$-presentation of $T$.
\end{proof}

\begin{remark} \label{rem-Q-E-type-res}
(i) In \cite{MDP-buch} Martin-Deschamps and Perrin have introduced
  $E$- and $N$-type resolutions of an ideal that are closely related
  to $E$-and $Q$-type resolutions as above. In fact,
$$
0 \to G_c \to \ldots \to G_2 \to Q \to I \to 0
$$
is an $N$-type resolution of the ideal $I$ if and only if
$$
0 \to G_c \to \ldots \to G_2 \to Q \to R \to R/I \to 0
$$
is a $Q$-type resolution of $R/I$. An analogous relation  is true for
the $E$-type resolutions of $I$ and $R/I$. In this sense, our
Definition \ref{def-E-and-Q-type} extends the concepts of $E$- and
$N$-type resolutions to modules with more than one generator.

(ii) As already indicated by the computation of cohomology modules
above, some properties of $M$ are
directly related to properties of the modules $E$ and $Q$,
respectively, in the corresponding resolutions of $M$.  For example, it
is easy to see that $E$ respectively $Q$ is a maximal
Cohen-Macaulay module if and only if $M$ is Cohen-Macaulay. If $M$ has
finite projective dimension then $M$ is Cohen-Macaulay if and only if
$E$ respectively $Q$ is a free module.

If $M$ is of pure codimension $c$ then $M$ is locally Cohen-Macaulay if
and only if it has cohomology of finite length and this is true if and
only if $E$ respectively $Q$ has cohomology of finite length. It follows
that in case $M$ has in addition finite projective dimension,
$M$ is (locally) Cohen-Macaulay if and only if $\tilde{E}$
respectively $\tilde{Q}$ is a vector bundle on $\Proj(R)$.
\end{remark}

A further relation between the modules $M, E, Q$ is stated in the
following result. It generalizes \cite{N-gorliaison}, Lemma 3.3.

Note that the module $E$ in an  $E$-type
resolution of an arbitrary module $M$ of codimension $c$ is always a
$c$-syzygy. If $M$ is unmixed  then it is even $(c+1)$-syzygy. More
precisely, we have.

\begin{lemma} \label{lem-unmixedness-by-E-and-Q-type-res}
Let $M$ be
  an $R$-module of codimension $c > 0$ having $E$- and
  $Q$-type resolution as in Definition \ref{def-E-and-Q-type}. Then
  the following conditions are equivalent:
\begin{itemize}
\item[(a)] $M$ is of pure codimension $c$.
\item[(b)] $Q$ is reflexive.
\item[(c)] $E$ is a $(c+1)$-syzygy.
\end{itemize}
\end{lemma}

\begin{proof}
Since reflexivity and being a $(c+1)$-syzygy can be cohomologically
 characterized
 (cf., e.g., \cite{N-gorliaison}, Proposition 2.5), our
claim follows by Lemma \ref{lem-char-unmixed} and the computation of
 cohomology in Remark
\ref{rem-to-def-of-E-and-Q-resolution} and Lemma
\ref{lem-Q-type-exists} .
\end{proof}

Now we are ready to show that resolutions of $E$- and $N$-type are
interchanged by direct m-linkage. The result generalizes Proposition 3.8 in
\cite{N-gorliaison}.

\begin{proposition} \label{prop-E-and-Q-resolution-under-linkage}
Let $M, N$ be $R$-modules of codimension $c > 0$ linked by the module
  $C$. Suppose $M$ has resolutions of $E$- and
  $Q$-type as in Definition \ref{def-E-and-Q-type}. Let
$$
0 \to D_c \to \ldots \to D_0 \to C \to 0
$$
be a minimal free resolution of $C$.  Put $s = r(C) + a(C) - r(R)$.   Then
$N$ has a $Q$-type resolution
$$
0 \to D'_c \oplus F_1^*(-s)  \to \ldots \to D_2 \oplus
F_{c-1}^*(-s)
\to D_1 \oplus E^*(-s) \to D_0 \to N \to 0
$$
where $D'_C$ is a free $R$-module such that $D'_c \oplus F_0^* \cong D_c$,
and an $E$-type resolution
$$
0 \to D''_c \oplus Q^*(-s) \to D_{c-1} \oplus G_2^*(-s) \to \ldots \to D_1 \oplus
G_c^*(-s) \to D_0 \to N \to 0
$$
where $D''_C$ is a free $R$-module such that $D''_c \oplus G_0^* \cong D_c$.
\end{proposition}

\begin{proof}
The proof is similar to the one of \cite{N-gorliaison}, Proposition
3.8. Thus we leave out some details which are treated there.
We proceed in several steps. We begin by showing the first claim
  starting with an $E$-type resolution of $M$ which we may and will assume to be minimal.

(I) Dualizing the given
$E$-type resolution of $M$ provides the complex
$$
0 \to R \to F_1^* \to \ldots \to F_{c-1}^* \to E^* \to \ER^c(R/I,R) \to 0
$$
which  is in fact an exact sequence.

Furthermore, we know by Lemma \ref{lem-Gor-versus-selfdual-res} that there are isomorphisms
$$
C \cong K_C (1 - r(C) - a(C)) \cong \ER^c(C, R) (-s).
$$
Thus, the self-duality of the minimal free resolution of $C$ means in particular that
$$
D_{c-i}^* \cong D_i (s) \quad \fall i = 0,\ldots,c.
$$

(II)
Lifting the homomorphism $\ffi: C \to M$ and using  Lemma \ref{lem-st-seq} we get a commutative diagram with exact rows and column
\begin{equation*}
\begin{CD}
&&&&&&&&&& 0 \\
&&&&&&&&&& @VVV \\
&&&&&&&&&& K_{N}(t) & \\
&&&&&&&&&& @VVV \\
0 @>>> D_c @>>> D_{c-1} @>>> \ldots @>>> D_0 @>>> C @>>> 0 \\
&& @VV{\ffi_c}V  @VV{\ffi_{c-1}}V && @VV{\ffi_0}V   @VV{\ffi}V \\
0 @>>> E @>>> F_{c-1} @>>> \ldots @>>> F_0 @>>> M @>>> 0. \\
&&&&&&&&&& @VVV \\
&&&&&&&&&& 0 \\
\end{CD}
\end{equation*}
Since the $E$-type resolution of $M$ is minimal, the homomorphism $\ffi_0$ is surjective. Thus, its $R$-dual $\ffi_0^* : F_0^* \to D_0^*$ is split-injective.

Now, dualizing the diagram above and using Step (I) we get by Definition \ref{def-linked-homo}
the commutative exact diagram
\begin{equation*}
\begin{CD}
&&&&&&&&&& 0 \\
&&&&&&&&&& @VVV \\
0 @>>> F_0^* @>>>  \ldots @>>> F_{c-1}^* @>>> E^* @>>> \ER^c(M, R) @>>> 0 \\
&& @VV{\psi}V  && @VVV @VVV   @VVV \\
0 @>>> D_c (s) @>>> \ldots @>>> D_1 (s) @>>> D_0 (s) @>>> C (s) @>>> 0  \\
&&&&&&&&&& @VV{L_C(\ffi)}V \\
&&&&&&&&&& N (s)  \\
&&&&&&&&&& @VVV \\
&&&&&&&&&& 0 \\
\end{CD}
\end{equation*}
where $\psi$ is the composition of $\ffi_0^*$ and an isomorphism. Hence, $\psi$ is split-injective, too. This shows that the module $F_0^*$ can be split off in the resulting mapping cone (cf.\ \cite{N-gorliaison}, Lemma 3.4). Thus, we get the exact sequence
$$
0 \to D'_c \oplus F_1^*(-s) \to  \ldots \to D_2 \oplus
F_{c-1}^*(-s)
\to D_1 \oplus E^*(-s) \to D_0 \to N \to 0.
$$
For it being a $Q$-type resolution, it remains to
show that $\HH^i(E^*) = 0$ if $n+2-c \leq i \leq n$.

 According to
Lemma \ref{lem-unmixedness-by-E-and-Q-type-res} we know that $E$ is a $(c+1)$-syzygy. Hence local duality and Lemma
\ref{lem-k-syzch} provide
$$
\HH^{n+1-i}(E^*)^{\vee}(1-r(R)) \cong \ER^i(E^*, R) = 0 \quad \mif 1 \leq i
\leq c-1.
$$
Thus, the argument for the $Q$-type resolution of $N$ is complete.

(III)
The proof for the $E$-type resolution of $N$ is similar. We only
sketch it. We may and will assume that the given $Q$-type resolution of $M$ is minimal.
 Replacing the $E$-type resolution of $M$ by the $Q$-type
resolution in the first diagram above and then dualizing provides the
following exact commutative diagram
$$
\xy\xymatrixrowsep{2.0pc}\xymatrixcolsep{1.4pc}\xymatrix{
&&&&&& 0 \ar @{->}[d]\\
0 \ar @{->}[r] & G_0^*  \ar @{->}[d]^-{\beta}  \ar @{->}[r] & Q^* \ar @{->}[d] \ar @{->}[r] & G_2^*  \ar @{->}[d] \ar @{->}[r] & \ldots \ar @{->}[r] & G_c^* \ar @{->}[d] \ar @{->}[r] & \ER^c(M, R)  \ar @{->}[d] \ar @{->}[r] & 0 \\
0 \ar @{->}[r] & D_c (s) \ar @{->}[r] & D_{c-1} (s) \ar @{->}[r] & D_{c-2}  \ar @{->}[r] & \ldots \ar @{->}[r] &  D_0 (s) \ar @{->}[r] & C (s) \ar @{->}[d]^-{L_c(\ffi)} \ar @{->}[r]  & 0  \\
&&&&&& N(s) \ar @{->}[d]\\
&&&&&& 0 
}
\endxy
$$
where $\beta$ is split-injective. Thus, we can split off $G_0^*$ in the mapping cone giving us
 the desired $E$-type resolution of $N$.
\end{proof}

In order to formulate some consequences of the last result we need more notation.

Let $M$ be an $R$-module of pure codimension $c \geq 1$. We have seen in Remark
\ref{rem-to-def-of-E-and-Q-resolution} and Lemma
\ref{lem-Q-type-exists} that the minimal  $E$- and
$N$-type resolution of  $M$ are uniquely determined.  Hence, there is a
well-defined map
$\ffi$ from the set of $R$-modules of pure codimension $c \geq 1$ into
the set of isomorphism classes of finitely generated $R$-modules
where $\ffi(M)$ is the class of the last module in a minimal $E$-type
resolution of $M$.

Similarly,   we get a well-defined map
$\psi$ from the
set of $R$-modules of pure codimension $c \geq 1$ into the set of
isomorphism classes of finitely generated $R$-modules  by defining
$\psi(M) = [Q]$ if $M$ has the
minimal $Q$-type resolution
$$
0 \to G_c \to \ldots \to  G_2 \to Q \to G_0  \to M \to 0.
$$

Recall that two  graded maximal $R$-modules $M$
  and $N$ are said to be {\itshape stably equivalent} if there are free
  $R$-modules $F, G$ and an integer $s$ such that
$$
M \oplus F \cong N(s) \oplus G.
$$
It is clear that stable equivalence is an equivalence relation.

Now we are able to state the main result of this section.

\begin{theorem}
  \label{even-liaison-class-versus-stable-euivalence-class} Let $c$ be
  a positive integer.
  The map $\ffi$ induces a well-defined  map $\Phi_c$ from the set
  $\cM_c$ of even
  liaison classes of modules of pure  codimension $c$ into the set
  $\cM_E^c$ of stable
  equivalence classes of finitely generated $(c+1)$-syzygies being
locally free in codimension $c-1$.

The map $\psi$ induces a well-defined map $\Psi_c$ from  $\cM_c$  into the set $\cM_Q^c$ of stable
  equivalence classes of finitely generated, reflexive modules $N$
that  satisfy $\HH^i(N) = 0$ for all $i$ with $n-c+2 \leq i \leq n$ and
  are locally free in codimension $c-1$.
\end{theorem}

\begin{proof} Proposition \ref{prop-E-and-Q-resolution-under-linkage}
  shows that
  the maps $\Phi_c$ and $\Psi_c$ do not depend on the choice of a
  representative of the even liaison class. If $M$ is a module of pure
  codimension $c$ then the localization of its $E$-type resolution at
  a prime $\fp \subset R$ of codimension $\leq c-1$ splits. Hence
  $\ffi(M)$ is locally free in codimension $c-1$. By Proposition
  \ref{prop-E-and-Q-resolution-under-linkage}, the same is true for
  $\psi(M)$. Thus,  Lemma \ref{lem-unmixedness-by-E-and-Q-type-res} and
Lemma \ref{lem-Q-type-exists} show that both maps $\Phi$ and $\Psi$
  are well defined
\end{proof}

The result above extends the analogous result for even Gorenstein
liaison classes of unmixed ideals (\cite{N-gorliaison}, Theorem 3.10) to
even module liaison classes.

\begin{remark} \label{rem-pol} If $R$ is just a polynomial ring over
  the field $K$ then the statement takes a somewhat simpler form
  because then every module in $\cM_Q^c$ and $\cM_E^c$ is
  automatically even locally free in codimension $c+1$. This follows
  from the fact that over a regular local ring $(c+1)$-syzygies are
  locally free in codimension $c+1$.
\end{remark}

\begin{remark} \label{comm-diagram-even-liason-stable-equivalence}
Using the notation in Theorem
\ref{even-liaison-class-versus-stable-euivalence-class} we have the
following commutative diagrams
\begin{equation*}
\begin{CD}
\cM_c @>{\Phi_c}>> \cM_E^c \\
@VV{\alpha}V @VV{\beta}V \\
\cM_c @>{\Psi_c}>> \cM_Q^c
\end{CD}
\end{equation*}
and
\begin{equation*}
\begin{CD}
\cM_c @>{\Psi_c}>> \cM_Q^c \\
@VV{\alpha}V @VV{\beta}V \\
\cM_c @>{\Phi_c}>> \cM_E^c
\end{CD}
\end{equation*}
where $\alpha$ is induced by linkage and $\beta$ is induced by dualization
with respect to $R$.
\end{remark}

Amasaki's main result in \cite{Amas-bourb} implies.

\begin{lemma} \label{lem-maps-surj}
If $R$ is a regular ring then the maps $\Phi$ and $\Psi$ in Theorem \ref{even-liaison-class-versus-stable-euivalence-class} are  surjective.
\end{lemma}

\begin{remark} \label{rem-maps-inj}
(i) The author expects that the preceding result is true without the assumption $R$ being regular. However, this requires new arguments because Amasaki's approach  heavily relies on the finiteness of free resolutions.

(ii) It remains a major challenge to decide whether the maps $\Phi$ and $\Psi$ are injective since an affirmative answer would provide a parametrization of even module liaison classes (cf.\ also Remark \ref{rem-mlicci})
\end{remark}

Theorem \ref{even-liaison-class-versus-stable-euivalence-class} implies,
for example, that in case $\ffi(M)$ and $\ffi(N)$ are not stably equivalent
the modules $M, N$ do not belong to the same even liaison
class. This shows that there is an abundance of even liaison classes,
but that there is also some control. This will be the topic of the following section.

We want to end this section by  discussing whether the module liaison
class of a given module $M$ contains a cyclic module. To this end we
recall that following Bruns (cf.\ \cite{Bruns-orientations} and
\cite{Bruns-alternating}), a finitely generated $R$-module $M$ is said
to be
{\it orientable} if it has a rank, is locally free in codimension one
and there is a homomorphism $\bigwedge^{\rank M} M \to R$ whose image has
codimension
at least two. Note that $M$ is orientable if it is locally free in codimension
one
and either $R$ is factorial or $M$ has finite projective dimension.

Theorem \ref{even-liaison-class-versus-stable-euivalence-class} has the following consequence.

\begin{corollary} \label{cor-whether-cyclic}
Let $M$ be a module of pure codimension $c \geq 2$. If there is a cyclic module in its even liaison class then  $M$ is orientable.
\end{corollary}

\begin{proof}
 This follows by the behavior of properties of orientable modules in exact sequences (\cite{Bruns-alternating}, Proposition 2.8). Indeed, if $N$ is a cyclic module then $\ffi(N)$ is orientable. Linking $N$ to another cyclic module we see that $\psi (N)$ is orientable, too. Now, Theorem \ref{even-liaison-class-versus-stable-euivalence-class} shows that all modules in the liaison class of $M$ are orientable.
\end{proof}

The last result raises the question whether  $M$ being orientable is not only a necessary, but also a sufficient condition for the liaison class of $M$ to contain a cyclic module.


\section{Transfer of  properties under liaison} \label{sec-transfer}

The goal of this section is to illustrate how the existence of the maps $\Phi$ and $\Psi$ can be used to show that cohomological and structural properties are preserved within (even) m-liaison classes. In particular, we generalize various results of Gorenstein liaison to our more general setting of module liaison.

We begin by discussing the local cohomology modules.

\begin{corollary} \label{cohomology_under_liaison} Let $M, N$ be modules of
  pure codimension $c$.
\begin{itemize}
\item[(a)] If $M$ and $N$ are in the same even liaison class then
  there is an integer $s$ such that
$$
\HH^i(M) \cong \HH^i(N) (s) \quad \fall i = 0,\ldots,n-c;
$$
\item[(b)] If $M$ is locally Cohen-Macaulay and if $M$ and $N$ are
  linked in an odd number of steps then there is an integer $s$ such that
$$
\HH^i(M) \cong \HH^{n+1-c-i}(N)^{\vee}(s) \quad \fall i = 1,\ldots,n-c;
$$
Moreover, if $M$ and $N$ are (directly) linked by the quasi-Gorenstein
module $C$ then $s = 1 -
r(C)- a(C)$.
\end{itemize}
\end{corollary}

\begin{proof} Part (a) is a consequence of  Theorem
  \ref{even-liaison-class-versus-stable-euivalence-class} and Remark
  \ref{rem-to-def-of-E-and-Q-resolution}.  It  remains to show the
  second claim of (b).  Let $E$ be a representative of the isomorphism class
  $\ffi(M)$. Then, using also Lemma \ref{lem-Q-type-exists}, we get
$$
\HH^i(E) \cong \HH^{i-c}(M) \quad \mif i \leq n
$$
and
$$
\HH^i(E^*)(r(R) - r(C) - a(C)) \cong \HH^{i-1}(N) \quad \mif i \leq n-c+1.
$$
Thus the claim is a consequence of local duality which provides
$$
\HH^i(E^*) \cong \HH^{n+2-i}(E)^{\vee}(1-r(R)) \quad \mif 2 \leq i
\leq n.
$$
\end{proof}

\begin{remark} \label{rem-coho-under-liai}
(i) The last result is an extension of the analogous result for
Gorenstein liaison classes of ideals(\cite{N-gorliaison}, Corollary 3.13).

(ii) Part (b) of the corollary above is not true if the modules are not locally Cohen-Macaulay. However, the intermediate cohomology modules of directly linked modules are related though  in general it seems difficult to make the relationship explicit. Chardin \cite{Chardin-dual} has some partial results in this direction for directly linked varieties of small dimension. These results can be extended to module linkage.
\end{remark}

Next, we consider the transfer of  structural properties  under module liaison.

\begin{corollary} \label{cor-CM-tranfer}
Let $M, N$ be  $R$-modules in the same module liaison class. Then we have:
\begin{itemize}
\item[(a)] $M$ is Cohen-Macaulay if and only if $N$ is Cohen-Macaulay.
\item[(b)] $M$ is \lCM\ if and only if $N$ has this property.
\end{itemize}
\end{corollary}

\begin{proof}
Claims (a) and (b) are immediate consequences of Corollary
\ref{cohomology_under_liaison} and the fact that $M$ is
Cohen-Macaulay, respectively locally Cohen-Macaulay if and only if the
cohomology modules $\HH^i(M),\ i < \dim M$, all vanish, respectively all have
finite length.
\end{proof}

A similar behavior is also true for Buchsbaum and surjective-Buchsbaum
modules. These
classes of modules properly
contain the class of Cohen-Macaulay modules, but cannot be
characterized by their local cohomology modules alone. For
comprehensive information about Buchsbaum modules, we refer to the
monograph \cite{St-V-buch} by St\"uckrad and
Vogel. Surjective-Buchsbaum modules have been introduced by
Yamagishi \cite{Yamagishi}. He observed that often Buchsbaum modules
are found by actually showing that they are even
surjective-Buchsbaum. Let us recall the definitions because we
use them later on.

Following Yamagishi \cite{Yamagishi}, the $R$-module
$M$ is called {\em surjective-Buchsbaum} if the natural homomorphisms
$\ffi^i_M : \Ext_R^i(K,M) \to \HH^i(M),\ i < \dim M$, are all surjective. Here
the maps $\ffi^i_M$ are induced by the embedding $0 :_M \fm \to
\HH^0(M)$.  Since $H^0(\fm, M) \cong 0 :_M \fm$ this embedding also
induces natural homomorphisms of derived functors $
\psi^i_M : H^i(\fm;M) \to \HH^i(M)
$ where $H^i(\fm ,M)$ is the $i$-th Koszul cohomology module of $M$
with respect to $\fm$. According to \cite{St-V-buch}, Theorem I.2.15,
the module $M$ is {\em Buchsbaum} if and only if $\psi^i_M$ is surjective
for all $i < \dim M$.

The isomorphism $H_0(\fm;R) = R/\fm \cong K$ lifts to a morphism of
complexes from the Koszul complex $K^{\bullet}(\fm;R)$ to a minimal
free resolution of $K$. It
induces natural homomorphisms
$
\lambda^i_M : \ER^i(K,M) \to H^i(\fm;M)
$.
Summing up, we have the following commutative diagram for all integers $i$
$$
\xy\xymatrixrowsep{0.5pc}\xymatrix{
\ER^i(K,M) \ar @{->}[dd]^-{\lambda^i_M}  \ar @{->}[dr]^-{\ffi^i_M}  \\
 & \HH^i(M). \\
H^i(\fm;M) \ar @{->}[ur]_-{\psi^i_M} }
\endxy
$$
The  diagram immediately shows that a
surjective-Buchs\-baum module is Buchsbaum. Note that the converse is
not true in
general. However, if $R$ is regular then $K_{\bullet}(\fm;M)$ is a minimal
free resolution of $K$, i.e., $\ER^i(K,M) \cong H^i(\fm;M)$. Hence, if $R$
is regular then an $R$-module is surjective-Buchsbaum if and only if it is
Buchsbaum.

The homological characterization of these modules allows us to trace their
properties along exact sequences. As a preparation, we need.

\begin{lemma} \label{lem-Buchs-E-Q-type}
Let $M$ be an $R$-module of codimension $c > 0$  and let $E, Q$ be
representatives of $\ffi (M)$ and $\psi(M)$, respectively. Then, if
one of the modules $M, E, Q$ is Buchsbaum or surjective-Buchsbaum then
all of them have the corresponding property.
\end{lemma}

\begin{proof}
We consider the Buchsbaum property first. Let
$$
0 \to T \to F \to M \to 0
$$
be an exact sequence of $R$-modules where $F$ is free. It induces the
following commutative diagram with exact rows
\begin{equation*}
\begin{CD}
H^i(\fm,F) @>>> H^i(\fm,M)  @>>> H^{i+1}(\fm,T)  @>>> H^{i+1}(\fm,F) \\
@VV{\psi^i_F}V @VV{\psi^i_M}V @VV{\psi^i_T}V @VV{\psi^{i+1}_F}V \\
\HH^i(F)  @>>>  \HH^i(M)  @>>> \HH^{i+1}(T)  @>>> \HH^{i+1}(F). \\
\end{CD}
\end{equation*}
Since the left-hand and the right-hand columns of this diagram vanish
if $i+1 < \dim R = n+1$ we get for every integer $k \geq 0$
that the map $\psi^i_M$ is surjective for all $i \leq k$ if and only
if $\psi^i_T$ is surjective for all $i \leq \min \{k+1, n\}$.

Consider now the $E$-type resolution of $M$
$$
0 \to E \to F_{c-1} \to \ldots \to F_0 \to M \to 0.
$$
Shopping it into short exact sequences the above observation shows
that $M$ is Buchsbaum if and only if $E$ is.

Next, consider the $Q$-type resolution of $M$
$$
0 \to G_c \to \ldots \to G_2 \to Q \to G_0 \to M \to 0
$$
where we may assume $c \geq 2$.
Reversing its construction in Lemma \ref{lem-Q-type-exists} we get the
exact sequences
$$
0 \to P \to Q \to T \to 0
$$
and
$$
0 \to T \to G_0 \to M \to 0
$$
where $P$ has projective dimension $c-2$, thus $\depth P = n+3-c$. The
first sequence induces the commutative diagram
\begin{equation*}
\begin{CD}
H^i(\fm,P) @>>> H^i(\fm,Q)  @>>> H^i(\fm,M)  @>>> H^{i+1}(\fm,P) \\
@VV{\psi^i_P}V @VV{\psi^i_Q}V @VV{\psi^i_M}V @VV{\psi^{i+1}_P}V \\
\HH^i(P)  @>>>  \HH^i(E)  @>>> \HH^i(M)  @>>> \HH^{i+1}(P). \\
\end{CD}
\end{equation*}
Using the vanishing of the cohomology of $Q$ in Lemma
\ref{lem-Q-type-exists} we always have that $\psi^i_Q$ is surjective
whenever $n+2 -c \leq i \leq n$. By the depth sensitivity of the
Koszul complex the left-hand and the right-hand columns of the diagram
vanish if $i \leq n+1-c$. We conclude that $Q$ is Buchsbaum if and
only if $\psi^i_T$ is surjective for all $i \leq n+1-c$ which, by the
first observation above, is equivalent to $M$ being Buchsbaum. This
completes the argument for the Buchsbaum property.

The proof for surjective-Buchsbaum modules is completely analogous. We
just have to replace the map $\psi^i_M$ by $\ffi^i_M$ everywhere in
the argument above.
\end{proof}

Now we want to use the Auslander dual in order to study Buchsbaumness
and surjective-Buchsbaumness under liaison. It allows us to simplify some arguments by avoiding the
use of derived categories.

The following result is essentially due to
St\"uckrad and Vogel.

\begin{lemma} \label{lem-A-dual}
Let $M$ be a maximal graded $R$-module with positive depth. Then:
\begin{itemize}
\item[(a)] If $M$ is a Buchsbaum module then $M^{\times}$ is
  so.
\item[(b)] If $M$ is a surjective-Buchsbaum module then  $M^{\times}$ is
  so.
\end{itemize}
\end{lemma}

\begin{proof}
Claim (a) is due to St\"uckrad and Vogel \cite{St-V-buch}, Proposition
III.1.28 as mentioned above. We sketch how the proof can be modified
to prove (b).

We may assume that $K$ is infinite. Then a sufficiently general linear
form $l \in R$ will be a non-zero divisor on $R, M$, and
$M^{\times}$. Set $\overline{M} := M/ l M$, $\overline{R} := R/ l R$ and denote
by $\overline{M}^{\times}$ the Auslander dual of $\overline{M}$ as
$\overline{R}$-module.

We will show the claim by induction on $n+1 = \dim M$. If $\dim M \leq
1$ then $M$, thus also $M^{\times}$ is Cohen-Macaulay. If $\dim M = 2$
then $M$ is surjective-Buchsbaum by (a) and \cite{N-loc-res}, Lemma
4.2, because $\depth M^{\times} > 0$.

Now let $\dim M \geq 3$. Then there is an isomorphism of $R$-modules
(cf.\ \cite{St-V-buch}, p.\ 173)
$$
\overline{M}^{\times} \cong (M^{\times}/ l M^{\times})/\HH^0 (M^{\times}/ l
    M^{\times})
$$
where $\HH^0 (M^{\times}/ l M^{\times})$ is annihilated by the maximal
  ideal $\fm$.
Since $M$ is surjective-Buchsbaum over $R$, $\overline{M}$ is
surjective-Buchsbaum over $\overline{R}$ by \cite{Yamagishi}, Theorem
3.2. Hence, by induction $\overline{M}^{\times}$ is a surjective-Buchsbaum
module over $\overline{R}$. Since $\fm \cdot  \HH^0 (M^{\times}/ l
  M^{\times}) = 0$, the isomorphism above implies that $M^{\times}/ l
    M^{\times}$ is a surjective-Buchsbaum module over $\overline{R}$. Using
    \cite{Yamagishi}, Theorem 3.2 again we conclude that $M^{\times}$
    is surjective-Buchsbaum over $R$.
\end{proof} 

We also need the following observation. 

\begin{lemma} \label{lem-stability-and-A-dual} 
Let $M$ and $N$ be directly linked maximal modules. If $M$ is not free then $N$ and $M^{\times}$ are stably equivalent. 
\end{lemma} 

\begin{proof} 
Let $M$ and $N$ be linked by the quasi-Gorenstein module $C$. Then $C$ must be free and there is an integer $t$ such that $C \cong C^* (t)$. Hence, there is   a minimal epimorphism $\pi: F \to M$, where $F$ is a free module, such that we get the following exact commutative diagram 
\begin{equation*}
\begin{CD}
&&&& 0 \\
&&&& @VVV \\ 
&&&& G \\
&&&& @VVV \\
0 @>>> \ker \ffi   @>>> C @>{\ffi}>> M  @>>> 0 \\
&& @VVV @VV{}V @VV{=}V \\
0 @>>> \ker \pi @>>> F  @>{\pi}>> M  @>>> 0 \\
&&&& @VVV \\ 
&&&& 0 
\end{CD}
\end{equation*}
where $G$ is a free module, too. Then, dualizing with respect to $R$
and shifting provide the exact commutative diagram 
\begin{equation*}
\begin{CD} 
&&&& 0\\
&&&& @VVV \\
0 @>>> M^*(t)  @>>> F^*(t)  @>>> M^{\times} (t) @>>> 0 \\
&& @VV{=}V @VV{}V @VVV \\
0 @>>> M^* (t) @>>> C  @>>> N  @>>> 0. \\
&&&& @VVV \\
&&&& G^* (t) \\
&&&& @VVV \\ 
&&&& 0
\end{CD}
\end{equation*}
Thus, the Snake lemma implies $N \cong M^{\times} (t) \oplus G^* (t)$,
completing the proof.  
\end{proof} 

Now we are ready to prove.

\begin{proposition} \label{prop-surj-Buchsb}
Let $M, N$ be modules in the same  liaison
  class. Then we have:
\begin{itemize}
\item[(a)] $M$ is Buchsbaum if and only if $N$ is so.
\item[(b)] $M$ is surjective-Buchsbaum
  if and only if $N$ is so.
\end{itemize}
\end{proposition}

\begin{proof} We may assume that $M$ and $N$ are directly linked by
  the quasi-Gorenstein module $C$. If $M$ is a free $R$-module then so is
  $N$. Thus, it suffices to consider non-free modules $M$ and $N$. 

(a) Suppose $M$ is Buchsbaum. We distinguish two cases.  First, assume
  that   $M$ is a maximal
  module. Then, by Lemma \ref{lem-stability-and-A-dual}, $N$ is stably
  equivalent to $M^{\times}$, thus Lemma \ref{lem-A-dual} gives the
  claim.

In this case, the linking module $C$ is a free
$R$-module. Thus, the exact 
sequence in Lemma \ref{lem-st-seq} shows that the $R$-dual $M^*$ of
$M$ is a Buchsbaum module, too. We will use this fact below.

Second, assume $\dim M < \dim R$. Let $E$ be a representative of $\ffi
(M)$ and let $Q$ be a representative of $\psi (N)$. Then Lemma
\ref{lem-Buchs-E-Q-type} shows that with $M$ also $E$ is Buchsbaum,
thus $E^*$ is Buchsbaum by the argument above. But Proposition
\ref{prop-E-and-Q-resolution-under-linkage} provides that $E^*$ and
$Q$ are stably equivalent. Hence, using Lemma \ref{lem-Buchs-E-Q-type}
again, we see that $N$ is Buchsbaum.

(b) By now it should be clear how this claim is proved analogously.
\end{proof}

\begin{remark} \label{rem-surj-Buchs} Part (a) of  Proposition
  \ref{prop-surj-Buchsb}  generalizes the corresponding result of
  Schenzel \cite{Schenzel-liai}  for Gorenstein liaison of ideals 
as well as the one of
Martsinkovsky and Strooker \cite{Marts-S} for their smaller module
  liaison classes. 
\end{remark}

Using $E$-type resolutions, Theorem \ref{even-liaison-class-versus-stable-euivalence-class} implies. 

\begin{lemma} \label{lem-proj-dim}
Let $M, N$ be modules in  the same even  liaison
  class. Then  $M$ has finite projective dimension if and only if $N$
  does.

Furthermore, $M$ and $N$ have the same projective dimension if it is finite.
\end{lemma}

Note that the analogous result is not true for the whole liaison class if $R$ is not
regular.

Abusing notation slightly, we say that $R/I$ is a {\em complete intersection} if $I$ is generated by an $R$-regular sequence. Note that every complete intersection is linked to itself by Example \ref{ex-mlink}(i). Thus,
Corollary \ref{cor-CM-tranfer} and Lemma \ref{lem-proj-dim} imply.

\begin{corollary} \label{cor-mlicci}
If $M$ is {\em mlicci} i.e.\ in the m-liaison class of a complete intersection, then $M$ is a perfect $R$-module.
\end{corollary}

\begin{remark} \label{rem-mlicci}
The converse of the last result would follow immediately if we knew
that the maps $\Phi$ and $\Psi$ in Theorem
\ref{even-liaison-class-versus-stable-euivalence-class} were
injective. However, we will show that the converse is true if the
codimension of the complete intersection is at most one (cf.\  Theorem
\ref{thm-codim1-perfect-mod}). 
\end{remark}

For modules of codimension zero, i.e.\ maximal modules, we can
describe their even liaison classes. 

\begin{proposition} \label{prop-l-class-of-max-mod} Let $M$ be an unmixed 
  maximal $R$-module. 
Then the module $N$ is in the even m-liaison class of $M$ if and only if 
$M$ and $N$ are stably equivalent. 
\end{proposition} 

\begin{proof} 
Let $N$ be a module in the even liaison class of $M$. We want to show that $M$ and $N$ are stably equivalent. This is clear if $M$ is free. Thus, we may assume that $M$ is not free and that $N$ is linked to $M$ in two steps. Let $P$ be a module that is directly linked to $M$ and $N$. Then, Lemma \ref{lem-stability-and-A-dual} shows that both $M$ and $N$ are stably equivalent to $P^{\times}$, hence $M$ and $N$ are stably equivalent, as claimed. 

For showing the reverse implication, let $N$ be a module that is stably equivalent to $M$. 
Applying Lemma \ref{lem-abspalten} with $D = R$ and also Lemma \ref{lem-shift} successively, we see that $N$ is in the even liaison class of $M$. 
\end{proof} 

Using Example \ref{ex-li-class-R}, we get in case $M = R$.  

\begin{corollary} \label{cor-free}
The module $N$ is in the m-liaison class of $R$ if and only if it is free.
\end{corollary}

In particular, over a field $K$ there is just one liaison class of $K$-modules.

\section{Liaison in codimension one} \label{sec-codim-1-liaison}

The goal of this section is to show that the perfect modules of codimension
one form the m-liaison class of the quotient ring of $R$ by a principal ideal.

\begin{theorem} \label{thm-codim1-perfect-mod}
Let $R$ be an integral domain and let $a \neq 0$ be an element of $R$ which
is not a unit.  Then an $R$-module $M$ belongs to the m-liaison
class of $R/a R$ if and only if $M$ is a perfect $R$-module of codimension
one.
\end{theorem}

Note that over an
integral domain a module is perfect of codimension one if and only if it
has a square presentation matrix with non-trivial determinant. Thus we will
deal with square matrices in the course of the proof.

We need some preparation and a bit of notation. Let $\ffi: F \to G$ be a
(graded) homomorphism between free modules represented by the homogeneous
matrix
$A$. Then we define $\coker A := \coker \ffi$.

The starting point is a special case of the result about the exchange of $E$-
and $Q$-type resolutions.

\begin{lemma} \label{lem-explicit-linked-module}
Let $F, G$ be (graded) free $R$-modules of the same rank and let
$\psi: G^*(s) \to F, \; 
\ffi:F \to G$ be (graded) homomorphisms which are not
isomorphisms. Choose bases for $F$ and $G$ and 
let $A , B$ be  the matrices representing $\ffi$ and $\psi$,
respectively. If $A \cdot B$ is equivalent to a (homogeneous)
symmetric matrix whose 
determinant is a non-zero divisor of $R$, then
$\coker \ffi$ and $\coker \psi^*(s)$ are m-linked by $\coker A B$.
\end{lemma}

\begin{proof} Put $S = A \cdot B$, $C = \coker S$ and $M = \coker
  A$. Since $\det A B = \det A \cdot \det B$ is a non-zero divisor of $R$
  there
  is a commutative diagram with exact rows
$$
\begin{array}{ccccccc}
  0 \longrightarrow & G^*(s) & \stackrel{S}{\longrightarrow} & G & \longrightarrow & C & \longrightarrow 0 \\
  & \Big\downarrow\! \! {\scriptstyle B} & & \| & & \Big\downarrow\! \! {\scriptstyle \gamma} & \\
0 \longrightarrow & F & \stackrel{A}{\longrightarrow} & G & \longrightarrow & M & \longrightarrow 0.
\end{array}
$$
Dualizing with respect to $R$ provides the exact commutative diagram
$$
\begin{array}{ccccccc}
  0 \longrightarrow & G^* & \stackrel{A^t}{\longrightarrow} & F^* & \longrightarrow & \ER^1(M, R) & \longrightarrow 0 \\
  & \| & & \Big\downarrow\! \! {\scriptstyle B^t} & & \Big\downarrow & \\
0 \longrightarrow & G^* & \stackrel{S^t}{\longrightarrow} & G(-s) & \longrightarrow & \ER^1(C, R) & \longrightarrow 0. \\
&&&&& \Big\downarrow \\
&&&&& \ER^1(\ker \gamma, R) \\
&&&&& \Big\downarrow \\
&&&&& 0
\end{array}
$$
By assumption, $S$ is equivalent to a symmetric matrix. Hence $C$ is a
quasi-Gorenstein module and $\ER^1(C,R)(s) \cong C$. Thus, we get $L_C(M)(-s) \cong \ER^1(\ker \gamma, R)$ by the definition of the
linking map. Therefore, the Snake lemma implies $\coker \psi^*(s) \cong L_C(M)$  completing the proof.
\end{proof}

This lemma suggests to introduce the notion of linked square matrices. Here
the restriction to Gorenstein rings is not necessary. Thus, we are working
in greater generality while dealing with matrices.

\begin{definition}
Let $R$ be an arbitrary ring. Then we denote the set of $n \times n$
matrices with entries in $R$ by $R^{n,n}$ and the transpose of a matrix $A$
by $A^t$. We say that two matrices $A, B
\in R^{n,n}$ are linked in one step if $A \cdot B^t$ is equivalent to a
symmetric matrix
whose determinant is a non-zero divisor of $R$. We call $A, B$ {\it linked
  matrices} if there are matrices $A = A_0, A_1,\ldots A_v = B$ such that
$A_i$ is linked in one step to $A_{i+1}$ for all $i = 0, 1,\ldots,v-1$. If
$R$ is a graded ring then we require additionally that all the matrices
$A_0,\ldots,A_v$ are homogeneous.
\end{definition}

It is obvious from the definition that being linked is an equivalence
relation among (homogeneous) square matrices of fixed size.

We will see that Theorem \ref{thm-codim1-perfect-mod} will essentially
follow from a result about linked matrices which we prove for more general
rings than Gorenstein rings.  Roughly speaking, the basic idea is to show
that over an integral domain a square matrix with non-vanishing
determinant is linked to a diagonal block matrix  with non-vanishing
determinant. In order to carry out this program we need two more preparatory
results.

\begin{lemma} \label{lem-preparation-of-matrix}
Let $R$ be an arbitrary integral domain. Furthermore,   in case $R =
\oplus_{i \geq 0} [R]_i$ is a graded ring assume that $[R]_1$ is
non-trivial. Let 
$A \in R^{n,n}\ (n \geq 2)$ be a square matrix  with non-vanishing
determinant
which is homogeneous if $R$ is graded.  Then there is a matrix
$\overline{A} := \left ( \begin{array}{cc}
a & \fb \\
\fc & A'
\end{array} \right ) \in R^{n, n}$ where $a \in R$ and $A' \in
R^{n-1,n-1}$ such that $\fb, 
\fc, \det \overline{A}$ are non-trivial and $\coker
A \cong \coker \overline{A}$. Furthermore, $\overline{A}$ can be taken
as a homogeneous matrix if $R$ is graded and $A$ is homogeneous. 
\end{lemma}

\begin{proof}
We have to show the existence of  invertible matrices $P, Q \in
R^{n,n}$ such that 
$\overline{A} = P A Q$ has the required properties. Performing suitable
elementary row and column operations on $A$, this is clear, at least if $R$ is not graded. It is a little more tricky if $R$ is graded because we have
less elementary row and column operations at our disposal. But, for example, an induction on $n$ will work. We omit the details. 
\end{proof}

\begin{lemma} \label{lem-existance-of-suitable-symmetric-matrices}
Let $R$ be a ring as in Lemma \ref{lem-preparation-of-matrix}.
Let $v, w \in R^n$ be non-trivial column vectors. Then there are a symmetric
matrix $S \in R^{n,n}$ and an element $\lambda \in R$ such that $\lambda
\neq 0, \det S \neq 0$ and $S v = \lambda w$. 

Furthermore, if $R$ is graded and $v = (v_1,\ldots,v_n)^t$, $w = (w_1,\ldots,w_n)^t$ are homogeneous such that $d := \deg v_i + \deg w_i$  for all $i = 1,\ldots,n$ then there are homogeneous $S$ and $\lambda$ with the properties above. 
\end{lemma}

\begin{proof} 
We restrict ourselves to the more difficult graded case. Then, by
assumption, $R$ contains  a linear form $L \neq 0$. Replacing all
powers of $L$ by the identity provides the argument in the non-graded
case.  

We begin with an observation which allows us to reduce the proof to the
most complicated case.

Suppose, for given vectors $v, w \in R^n$ we have found $\lambda$ and
$S$ as in the 
statement. Consider the vectors
$$
v' = \left(\begin{array}{c}
v_0 \\
v
\end{array} \right ), \; w' = \left(\begin{array}{c}
w_0 \\
w
\end{array} \right ) \in R^{n+1}.
$$
In case that both $v_0$ and $w_0$ are non-trivial, we get the desired
conclusion for $v', w'$  because putting 
$$
S' = \left (\begin{array}{cc}
\lambda w_0 & 0 \\
0 & S v_0
\end{array} \right ) \in R^{n+1, n+1}
$$
we obtain
$$
S' v' = (\lambda v_0) w'
$$
where $\det S', \ \lambda v_0 \neq 0$. 

Assume now that we have $v_0 = w_0 = 0$.  Multiplication by $S$
induces a homomorphism $G \to G^*(s)$ where $G$ is a graded free
$R$-module of rank $n$ and $s \in \ZZ$. Since $v_0 = w_0 = 0$ we may choose
$d_0 := \deg v_0$ such that $s - 2 d_0 \in \{0, 1\}$.  Then the
conclusion of the statement follows for $v', w'$ because $S' v' = 
\lambda w'$ where $S'$ is the homogenous matrix
$$
S' = \left (\begin{array}{cc}
L^{s - 2 d_0} & 0 \\
0 & S
\end{array} \right ) \in R^{n+1, n+1}.
$$

Using the observation above (and possibly reordering the rows)
we see that it suffices to show the statement for vectors
$$
v = (0,\ldots, 0, v_{k+1},\ldots, v_n)^t, \; w = (w_1,\ldots, w_k,
0,\ldots, 0)^t 
$$
where $k$ is an integer with $1 \leq k < n$  and all entries
$v_{k+1},\ldots, v_n, w_1,\ldots, w_k$ are non-trivial. In this
situation, we can always adjust the degrees of the entries of $v, w$
such that the degree assumption is satisfied and, in particular, we can
choose $d$ sufficiently large. 

Now we distinguish two cases. 

\noindent
{\it Case 1}. Assume $k \geq \frac{n}{2}$.

Put $\lambda = v_{k+1} \cdot
\ldots \cdot v_n$. The corresponding product where one factor $v_j$ is
omitted will 
be abbreviated by $\frac{\lambda}{v_j} \in R$. Consider the following matrices
$$
A = \left ( \begin{array}{ccccc}
\frac{\lambda}{v_{k+1}} w_1 & & & &\\
& \ddots &   & 0 & \\
&&&& \\
&&& \frac{\lambda}{v_{n-1}} w_{n-k-1} & \\[5pt]
&&&& \frac{\lambda}{v_n} w_{n-k} \\[5pt]
&&&& \frac{\lambda}{v_n} w_{n-k+1} \\
&& 0 && \vdots \\
&&&& \frac{\lambda}{v_n} w_k
\end{array} \right ) \in R^{k, n-k}
$$
and
$$
S = \left ( \begin{array}{ccc|c}
0 & & 0 & \\
& & & A \\
0 & & D & \\[3pt] \hline
& & & \\[-1.5ex]
& A^t && 0
\end{array} \right ) \in R^{n,n}
$$
where $D$ denotes the diagonal $(2k-n) \times (2k-n)$ matrix whose
$j$-th entry on the main diagonal is $L$ to the power $d + \deg
\lambda - 2 \deg v_{n-k+j}$. Here, we chosse $d$ large enough such
that all the powers of $L$ have a non-negative exponent. 
It is easy to check
that $S$ is a homogeneous matrix,  
$$
S v = \lambda w, 
$$
and
$$
\det S = \pm \left ( \prod_{i=1}^{n-k} \frac{\lambda}{v_{k+i}} w_i \right )
\cdot
  \det \left (
\begin{array}{ccc}
0 & & D \\[3pt] \hline \\[-1.5ex]
& A^t &
\end{array} \right ) = \pm \left ( \prod_{i=1}^{n-k} \frac{\lambda}{v_{k+i}}
w_i \right )^2 \cdot L^e
\neq 0
$$
for some $e \in \ZZ$,  whence the claim.

\noindent
{\it Case 2}. Assume $k \leq \frac{n}{2}$.

Applying Case 1 we find a matrix $S$ and $\lambda \in R$ such that $\det S, \
\lambda \neq 0$ and $S w = \lambda v$. Multiplying the  last equation by
the adjoint matrix of $S$ we obtain
$$
\det S \cdot w = \lambda \cdot  \adj S \cdot v
$$
which proves the claim because $\adj S$ is symmetric if $S$ is a symmetric
matrix.
\end{proof}

Now we are ready for the announced result about linked matrices.

\begin{lemma} \label{lem-key-for-linked-matrices}
Let $R$ be a ring as in Lemma \ref{lem-preparation-of-matrix}. Let $A =
\left ( \begin{array}{cc} 
a & \fb \\
\fc & A'
\end{array} \right ) \in R^{n,n}$ be a
square matrix where $a \in R,\; \fc, \fb^t \in R^{n-1},\; A' \in R^{n-1, n-1}$. 
If  $\det A, \det A', \fb$, and 
$\fc$ are non-trivial then $A$ is linked to a  square matrix $\left ( \begin{array}{cc}
b & 0 \\
0 & B'
\end{array} \right )$. 

Furthermore, $\left ( \begin{array}{cc}
b & 0 \\
0 & B'
\end{array} \right )$ can be taken as a homogeneous matrix if $R$ is
graded and $A$ is homogeneous. 
\end{lemma}

\begin{proof}
Put $\tilde{\fb} = \fb \cdot \adj A'$ where $\adj A'$ denotes the adjoint
matrix of $A'$. Then $\tilde{\fb}$ is non-trivial because otherwise we
would get
$$
0 = \tilde{\fb} \cdot A' = \fb \cdot \adj A' \cdot A = \fb \cdot \det A'
$$
which is a contradiction since $\fb$ and $\det A'$ are non-trivial by
assumption.

Thus we can apply Lemma \ref{lem-existance-of-suitable-symmetric-matrices}
and conclude that there are a symmetric
matrix $\tilde{S} \in R^{n,n}$ and an element $\lambda \in R$ such that
$\lambda
\neq 0, \det \tilde{S} \neq 0$ and $\tilde{\fb} \tilde{S}  = \lambda
\fc^t$.

Now we define the matrices $B \in R^{n,n}$ and $B' \in R^{n-1,n-1}$  by
$$
B' := \adj A' \cdot \tilde{S} \quad \mbox{and} \quad
B^t := \left ( \begin{array}{cc}
\lambda & 0 \\
0 & B'
\end{array} \right ).
$$
It follows that
$$
S := A \cdot B^t = \left ( \begin{array}{cc}
a \lambda & \fb \cdot B' \\
\lambda \fc  & A' \cdot B'
\end{array} \right )
$$
which is a symmetric matrix because
$$
A' \cdot B' = A' \cdot \adj A' \cdot \tilde{S} = \det A' \cdot \tilde{S}
$$
is symmetric and
$$
\lambda \fc^t = \tilde{\fb} \tilde{S} = \fb \cdot \adj A' \cdot \tilde{S} =
\fb \cdot B'
$$
due to our choice of $\tilde{S}$. Furthermore, $S$ has non-trivial
determinant since
$$
\det S = \det A \cdot \lambda \cdot \det(\adj A') \cdot \det \tilde{S}
$$
and each factor on the right-hand side is non-trivial. Therefore, the
matrices $A$ and $B$ are linked and we are done.
\end{proof}

Now we are in a position to show the main result of this section.

\begin{proof}[Proof of Theorem \ref{thm-codim1-perfect-mod}]
 One direction is clear by Corollary \ref{cor-mlicci}.

In order to show the converse, let $A \in R^{n,n}$ be a presentation matrix of
$M$. If $n = 1$ there is nothing to show. Let $n \geq
  2$. According to Lemma \ref{lem-preparation-of-matrix} we may assume that
  $A = \left ( \begin{array}{cc}
a & \fb \\
\fc & A'
\end{array} \right )$ has the property that $\fb, \fc$ and $\det A'$ are
non-trivial. Lemma \ref{lem-key-for-linked-matrices} shows that there is a
matrix $B = \left ( \begin{array}{cc}
b & 0 \\
0 & B'
\end{array} \right )$ which is linked to $A$. In spite of Lemma
\ref{lem-explicit-linked-module} 
we obtain that the modules $M$ and $\coker B$ are linked. By Lemma
\ref{lem-abspalten}, 
it follows that $\coker B$ and $\coker B'$ are evenly
linked. Altogether we obtain that 
$M = \coker A$ is in the same m-liaison class as $\coker B'$. Thus we
conclude by induction on $n$ that $M$ is in the m-liaison class of
$(R/c R)(j)$ for some $j \in \ZZ$ and 
some $c \neq 0$. The module $(R/c R)(j)$ is linked to $(R/a R)(j)$ by $(R/a c
R)(j)$. Now, $(R/a R)(j)$ and $R/a R$ are in the same even liaison
class by Lemma 
\ref{lem-shift}. This completes the argument.
\end{proof}

\end{document}